\documentclass[12pt]{amsart}

\usepackage{amssymb,amsmath}
\usepackage[ps2pdf,colorlinks=true,urlcolor=blue,
citecolor=red,linkcolor=blue,linktocpage,pdfpagelabels,bookmarksnumbered,bookmarksopen]{hyperref}
\usepackage{epsfig,graphics}

\usepackage[left=2.6cm,right=2.6cm,top=2.6cm,bottom=2.6cm]{geometry}

\newcommand{\C}{{\mathbb C}}
\newcommand{\R}{{\mathbb R}}

\newcommand{\RN}{{\mathbb{R}^{N}}}
\newcommand{\rn}{{\mathbb{R}^N}}
\newcommand{\be}{\begin{equation}}
\newcommand{\ee}{\end{equation}}

\renewcommand{\H}{{\mathbb H}}
\newcommand{\elle}[1]{L^{#1}}

\renewcommand{\b }{\beta }

\newcommand{\teta }{\theta }
\newcommand{\vep}{\varepsilon}
\newcommand{\eps}{\varepsilon}

\newcommand{\E}{{\mathcal E}}
\newcommand{\calo}{{\mathcal O}}
\newcommand{\cala}{{\mathcal A}}

\newcommand{\im}{{\mathcal Im}}
\newcommand{\re}{{\mathcal Re}}
\numberwithin{equation}{section}
\newtheorem{theorem}{Theorem}[section]
\newtheorem{proposition}[theorem]{Proposition}

\newtheorem{lemma}[theorem]{Lemma}
\newtheorem{definition}[theorem]{Definition}
\theoremstyle{definition}
\newtheorem{remark}[theorem]{Remark}

\newcommand{\brm}{\begin{remark}\rm}
\newcommand{\erm}{\end{remark}}
\newcommand{\brms}{\begin{remark}\rm}
\newcommand{\erms}{\end{remark}}
\newcommand{\bte}{\begin{theorem}}
\newcommand{\ete}{\end{theorem}}
\newcommand{\bpr}{\begin{proposition}}
\newcommand{\epr}{\end{proposition}}
\newcommand{\ble}{\begin{lemma}}
\newcommand{\ele}{\end{lemma}}
\newcommand{\beq}{\begin{equation}}
\newcommand{\eeq}{\end{equation}}
\newcommand{\bdm}{\begin{displaymath}}
\newcommand{\edm}{\end{displaymath}}
\numberwithin{equation}{section}

\newcommand{\bos}{\begin{remark}\rm}
\newcommand{\eos}{\end{remark}}
\newcommand{\dys}{\displaystyle}
\newcommand{\lbl}{\label}
\newcommand{\bdim}{{\noindent{\bf Proof.}\,\,\,}}
\newcommand{\edim}{{\unskip\nobreak\hfil\penalty50\hskip2em\hbox{}\nobreak\hfil\mbox{\rule{1ex}{1ex} \qquad}
  \parfillskip=0pt \finalhyphendemerits=0\par\medskip}}
\newcommand{\ben}{\begin{enumerate}}
\newcommand{\een}{\end{enumerate}}

\title[Soliton dynamics for CNLS systems with potentials]
{Soliton dynamics for CNLS systems with potentials}

\author[E.\ Montefusco]{Eugenio Montefusco}
\author[B.\ Pellacci]{Benedetta Pellacci}
\author[M.\ Squassina]{Marco Squassina }

\address{Dipartimento di Matematica
\newline\indent
Sa\-pien\-za Universit\`a di Roma
\newline\indent
Piazzale A.\ Moro 5, I-00185 Roma, Italy}
\email{montefusco@mat.uniroma1.it}

\address{Dipartimento di Scienze Applicate
\newline\indent
Universit\`a degli Studi di Napoli ``Parthenope''
\newline\indent CDN Isola C4, I-80143 Napoli, Italy}
\email{benedetta.pellacci@uniparthenope.it}

\address{Dipartimento di Informatica
\newline\indent
Universit\`a degli Studi di Verona
\newline\indent
C\'a Vignal 2, Strada Le Grazie 15, I-37134 Verona, Italy}
\email{marco.squassina@univr.it}

\thanks{The first and the second author are supported by the MIUR national
research project ``Variational Methods and Nonlinear Differential Equations'',
while the third author is supported by the 2007 MIUR national research
project ``Variational and Topological Methods in the Study of
Nonlinear Phenomena''}
\begin{document}

\begin{abstract}
The semiclassical limit of a weakly coupled  nonlinear focusing Schr\"odinger system
in presence of a nonconstant potential is studied.
The initial data is of the form $(u_{1},u_{2})$ with $u_{i}=r_{i}\big(
\frac{x-\tilde x}{\vep}\big)e^{\frac{{\rm i}}{\vep} x\cdot \tilde{\xi}}$, 
where $(r_{1},r_{2}) $  is a real ground state solution, belonging to a suitable class, of an associated autonomous 
elliptic system. For $\vep$ sufficiently small, 
the solution $(\phi_{1},\phi_{2})$ will been shown to have, locally in time, the form  
$(r_{1}\big(\frac{x- x(t)}{\vep}\big) e^{\frac{{\rm i}}{\vep} 
x\cdot {\xi}(t)},r_{2}\big(\frac{x- x(t)}{\vep}\big)
e^{\frac{{\rm i}}{\vep} x\cdot {\xi}(t)})$,
where $(x(t),\xi(t))$ is the solution of 
the Hamiltonian system $\dot x(t)=\xi(t)$, $\dot \xi(t)=-\nabla V(x(t))$
with $x(0)=\tilde{x}$ and $\xi(0)=\tilde{\xi}$.
\end{abstract}
\maketitle


\section{Introduction and main result}
 
\subsection{Introduction}
In recent years much interest has been devoted to the study of systems of weakly coupled nonlinear  Schr\"odinger equations. This interest is motivated 
by many physical experiments especially in nonlinear optics and in the theory of Bose-Einstein condensates
(see e.g.\ \cite{a,hk,man,men}). Existence results of {\em ground and bound states} solutions have 
been obtained by different authors (see e.g.\ \cite{aCol,BartWa,defiguerlopes,lw,mmp1,sirakov}).
A very interesting  aspect regards the dynamics, in the semiclassical limit, of a {\em general solution},
that is to consider  the nonlinear  Schr\"odinger system 
\beq
\label{prob01}
\begin{cases}
\dys i\vep \partial_{t}\phi_{1}+\frac{\vep^{2}}{2}\Delta \phi_{1}
- V(x)\phi_{1}+\phi_{1}(|\phi_{1}|^{2p}
+\b|\phi_{2}|^{p+1}|\phi_{1}|^{p-1})=0 
& \text{in $\R^N\times\R^+$}, 
\\
\noalign{\vskip4pt}
\dys i\vep \partial_{t}\phi_{2}+\frac{\eps^2}2\Delta \phi_{2}
 - V(x)\phi_{2}+\phi_{2}(|\phi_{2}|^{2p}
+\b|\phi_{1}|^{p+1}|\phi_{2}|^{p-1})=0
 & \text{in $\R^N\times\R^+$},
\\
\noalign{\vskip4pt}
\dys \phi_{1}(0,x)=\phi_{1}^{0}(x)\qquad
\dys \phi_{2}(0,x)=\phi_{2}^{0}(x),
\end{cases}
\end{equation}
with $0<p<2/N$, $N\geq 1$ and $\beta>0$ is a constant modeling the birefringence effect of the material. The potential $V(x)$ is a regular function in $\R^{N}$  modeling the action of external forces (see \eqref{ipoV}),
$\phi_{i}:\R^+\times \RN\to \C$ are complex valued
functions and $\vep>0$ is a small parameter playing the r\^ole of Planck's constant. 
The task to be tackled with respect to this system is to recover the full dynamics of
a solution $(\phi_{1}^{\vep},\phi_{2}^{\vep})$ as a point particle subjected
to galileian motion for the parameter $\vep$ sufficiently small.
Since the famous papers~\cite{abc,dpf,fw}, a large amount of work has been dedicated to this study 
in the case of a single Schr\"odinger equation and for a special class of solutions, 
namely standing wave solutions (see~\cite{ambook} and the references therein). When considering 
this particular kind of solutions one is naturally lead to study the following elliptic system 
corresponding to the physically relevant case $p=1$ (that is Kerr nonlinearities)
\begin{equation}
\label{kerrsystem}
\begin{cases}
-\eps^2\Delta u + V(x)u =u^3+\beta v^2u & \text{in $\R^N$}, \\
\noalign{\vskip2pt}
-\eps^2\Delta v + V(x)v =v^3+\beta u^2v & \text{in $\R^N$},
\end{cases}
\end{equation}
so that the analysis reduces to the study of the asymptotic behavior of solutions of an elliptic system. The concentration of a  least energy  solution  around the local minima (possibly degenerate) of the potential $V$ has been studied in~\cite{mps}, where some sufficient and necessary conditions  have been
established. To our knowledge the semiclassical dynamics
of different kinds of solutions of a single Schr\"odinger equation 
has been tackled in the series of papers~\cite{bronski, Keerani1, Keerani2} (see also
\cite{BeghMi} for recent developments on the long term soliton dynamics), 
assuming that the initial datum is of the form $r((x-\tilde{x})/\vep)e^{{\frac{\rm i}{\vep}}x\cdot \tilde\xi}$, 
where $r$ is the unique ground state solution of 
an associated elliptic problem (see equation~\eqref{scalareqq}) and $\tilde x,\tilde\xi\in\R^N$. 
This choice of initial data corresponds to the study of a different situation from the previous one. Indeed, 
it is taken into consideration the semiclassical dynamics of ground state solutions of the autonomous 
elliptic equation once the action of external forces occurs.
In these papers it is proved that the solution is approximated by the ground state $r$--up to translations and phase changes--and the translations  and phase changes are precisely related with the solution of a Newtonian system in $\R^N$ governed by the 
gradient of the potential $V$.
Here we want to recover similar results for system~\eqref{prob01} taking as initial data
\beq
\lbl{dati}
\phi_{1}^{0}(x)=r_{1}\Big(\frac{x-\tilde x}{\vep}\Big)e^{\frac{{\rm i}}{\vep} x\cdot \tilde{\xi}},\qquad 
\phi_{2}^{0}(x)=r_{2}\Big(\frac{x-\tilde x}{\vep}\Big)e^{\frac{{\rm i}}{\vep} x\cdot \tilde{\xi}},
\eeq
where the vector $R=(r_{1},r_{2})$ 
is a suitable ground state (see Definition~\ref{admissible}) of the associated elliptic system
\beq
\label{probel0}
\tag{$E$}
\begin{cases}
\dys -\frac{1}{2}\Delta r_{1}+r_{1}=r_{1}(|r_{1}|^{2p}+
\b|r_{2}|^{p+1}|r_{1}|^{p-1})
& \text{in $\R^N$}, 
\\
\noalign{\vskip4pt}
\dys -\frac{1}2\Delta r_{2} +r_{2}=r_{2}(|r_{2}|^{2p}
+\b|r_{1}|^{p+1}|r_{2}|^{p-1})
 & \text{in $\R^N$}.
\end{cases}
\end{equation}
When studying the dynamics of systems some new difficulties can arise. 
First of all, we have to take into account that, up to now, it is still not known if a uniqueness result  (up to
translations in $\R^N$) for real ground state solutions of~\eqref{probel0} holds. This is expected, at least in the case where $\beta>1$. Besides, also nondegeneracy properties (in the sense provided
in~\cite{dw,mps-modstabil}) are proved in some particular cases~\cite{dw,mps-modstabil}.
These obstacles lead us to restrict the set of admissible ground state solutions we will take into consideration (see Definition~\ref{admissible}) in the study of soliton dynamics.

Our first main result (Theorem~\ref{mainth3}) will give the desired asymptotic behaviour. 
Indeed, we will show that a solution which starts from~\eqref{dati} (for a suitable ground state $R$) will
remain close to the set of ground state solutions, up to translations and phase rotations. 
Furthermore, in the second result (Theorem~\ref{mainth-Bis}), we will prove that the mass densities 
associated with the solution $\phi_i$ converge--in the dual space of $C^{2}(\rn)\times C^{2}(\rn)$--to the delta measure with  mass given by $\|r_{i}\|_{L^{2}}$ and concentrated along $x(t)$, solution to the (driving) Newtonian differential equation  
\beq\lbl{Ham1}
\ddot{x}(t)=-\nabla V(x(t)),\qquad x(0)=\tilde{x},\;\;
\dot{x}(0)=\tilde{\xi}
\eeq
where $\tilde{x}$ and $\tilde{\xi}$ are fixed in the initial data of \eqref{prob01}.
A similar result for  each single component of the momentum density is lost as a consequence of the birefringence effect. However, we can afford the desired result for a balance on the total momentum density. 
This shows that--in the semiclassical regime--the solution moves as a point particle 
under the galileian law given by the Hamiltonian system~\eqref{Ham1}.
In the case of $V$ constant our statements are related with the results 
obtained, by linearization procedure, in \cite{weinsteinMS} for the single 
equation. Here, by a different approach, we show that \eqref{Ham1} gives
a modulation equation for the solution generated by the initial data
\eqref{dati}. Although we cannot predict the shape of the solution,
we know that the dynamic of the mass center is described by \eqref{Ham1}.
The arguments will follow~\cite{bronski, Keerani1, Keerani2}, where the case of a single Schr\"odinger equations has been considered. 
The main ingredients are the {\em conservation laws} of~\eqref{prob01} and 
of the Hamiltonian associated with the ODE in~\eqref{Ham1} and a 
{\em modulational stability} property for a suitable class of ground state solutions for the associated 
autonomous elliptic system~\eqref{probel0}, recently proved in~\cite{mps-modstabil} by the authors
in the same spirit of the works~\cite{weinsteinMS,weinsteinCpam}
on scalar Schr\"odinger equations.

The problem for the single equation has been also studied using the 
WKB analysis (see for example~\cite{c} and the references therein), to our knowledge, there are no results for the system using this approach.
Some of the arguments and estimates in the paper are strongly based 
upon those of~\cite{Keerani2}. On the other hand,  for the sake of
self-containedness, we prefer to include all the details in the proofs.  
\smallskip

\subsection{Admissible ground state solutions}

Let $\H_{\vep}$ be the space of the vectors $\Phi=(\phi_{1},\phi_{2})$ 
in $\H=H^1(\R^N;\C^2)$ endowed with the rescaled norm
$$
\|\Phi\|^{2}_{\H_{\vep}}=
\frac1{\vep^{N}}\|\Phi\|_{2}^{2}+
\frac1{\vep^{N-2}}\|\nabla \Phi\|_{2}^{2},
$$
where 
$\|\Phi\|_{2}^{2}=\|(\phi_{1},\phi_{2})\|_{2}^{2}=
\|\phi_{1}\|_{2}^{2}+\|\phi_{2}\|_{2}^{2}$ 
and $\|\phi_{i}\|_{2}^{2}=\|\phi_i\|_{\elle2}^{2}$ is the standard norm in the Lebesgue space $L^{2}$ given by $\|\phi_i\|_{2}^{2}=\int\phi_i(x)\bar\phi_i(x)dx$.
\\
We aim to study the semiclassical dynamics of a least energy solution of problem $\eqref{probel0}$ once the action of external forces is taken into consideration. 

In~\cite{aCol,mmp1,sirakov} it is proved that there exists a least action
solution $R=(r_1,r_2)\neq(0,0)$ of \eqref{probel0} which has 
nonnegative components. Moreover, $R$ is a solution to
the following minimization problem (cf.\ \cite[Theorems 3.4 and 3.6]{mmp2-pp})
\beq\label{caravar}
{\mathcal E}(R)=\min_{{\mathcal M}}{\mathcal E},\quad\text{where}\quad
{\mathcal M}=\left\{ U\in {\H} : \|U\|_{2}=\|R\|_{2}\right\}, 
\eeq
where the functional $\E:\H\to \R$ is defined by 
\begin{align}\label{Ecorsivo}
\E(U)&=\frac12\|\nabla U\|_{2}^{2}-\int F_{\beta}(U)dx
\\
\label{Fbeta} 
F_{\beta}(U)&=\frac1{p+1}\Big(|u_{1}|^{2p+2}+|u_{2}|^{2p+2}+2\beta
|u_{1}|^{p+1}|u_{2}|^{p+1}\Big),
\end{align}
for any $U=(u_1,u_2)\in\H$.
We shall denote with ${\mathcal G}$ the set of the (complex) ground 
state solutions.

\begin{remark}
\label{RemRaprGS}
Any element $V=(v_1,v_2)$ of $\mathcal{G}$ has the form
\begin{equation*}
V(x)=(e^{i\theta_1}|v_1(x)|,e^{i\theta_2}|v_2(x)|),\quad x\in\R^N,
\end{equation*}
for some $\theta_1,\theta_2\in S^{1}$ (so that $(|v_1|,|v_2|)$ is a real, positive, ground state solution). 
Indeed, if we consider the minimization problems
\begin{align*}
\sigma_\C&=\inf\big\{\E(V): V\in \H,
\,\|V\|_{\elle2}=\|R\|_{\elle2} \big\},\\
\sigma_\R&=\inf\big\{\E(V): V\in H^1(\R^N;\R^{2})\,
\|V\|_{\elle2}=\|R\|_{\elle2} \big\}
\end{align*}
it results that $\sigma_\C=\sigma_\R$.
Trivially one has $\sigma_\C\leq \sigma_\R$. Moreover, if $V=(v_1,v_2)\in 
\H$, due to the well-known pointwise inequality $|\nabla |v_i(x)||\leq |\nabla v_i(x)|$ 
for a.e.\ $x\in\R^N$, it holds
$$
\int|\nabla |v_i(x)||^2dx\leq \int|\nabla v_i(x)|^2dx,\quad i=1,2,
$$
so that also $\E(|v_1|,|v_2|)\leq \E(V)$. In particular, we conclude that 
$\sigma_\R\leq\sigma_\C$,
yielding the desired equality 
$\sigma_\C=\sigma_\R$. Let now $V=(v_1,v_2)$ be a {\em 
solution} to $\sigma_\C$
and assume by contradiction that, for some $i=1,2$,
$$
{\mathcal L}^N(\{x\in\R^N:|\nabla |v_i|(x)|<|\nabla v_i(x)|\})>0,
$$
where ${\mathcal L}^N$ is the Lebesgue measure in $\R^N$.
Then $\|(|v_1|,|v_2|)\|_{L^2}=\|V\|_{L^2}$, and
\begin{equation*}
\sigma_\R \leq \frac{1}{2}\sum_{i=1}^2\int |\nabla |v_i||^2dx-\int F_\beta(|v_1|,|v_2|)dx
<\frac{1}{2}\sum_{i=1}^2\int |\nabla v_i|^2dx-\int F_\beta(v_1,v_2)dx=\sigma_\C,
\end{equation*}
which is a contradiction, being $\sigma_\C=\sigma_\R$. Hence, we have
$|\nabla |v_i(x)||=|\nabla v_i(x)|$ for a.e.\ $x\in\R^N$ and any $i=1,2$. This is true
if and only if $\re\, v_i\nabla (\im\, v_i)=\im\, v_i\nabla(\re\, v_i)$.
In turn, if this last condition holds, we get
$$
{\bar v_i}\nabla v_i=\re\, v_i\nabla (\re\, v_i)+
\im\, v_i\nabla (\im\, v_i),\quad \text{a.e.\ in $\R^N$},
$$
which implies that $\re\,(i\bar v_i(x)\nabla v_i(x))=0$ a.e.\ in $\R^N$. Finally, for any 
$i=1,2$, from this last identity one immediately
finds $\theta_i\in S^1$ with $v_i=e^{i\theta_i}|v_i|$, concluding the proof.
\end{remark}

	In the scalar case, the ground state solution for the equation
	\begin{equation}
		\label{scalareqq}
	-\frac{1}{2}\Delta r+r=r^{2p+1}\quad\text{in $\R^N$}
\end{equation}
	is always unique (up to translations) and nondegenerate (see e.g.~\cite{kwong,MlSe,weinsteinMS}).
	For system~\eqref{probel0}, in general, the uniqueness and nondegeneracy of ground state
	solutions is a delicate open question. 
	
	The so called {\em modulational stability} property of ground states solutions 
plays an important r\^ole in soliton dynamics on finite time intervals.
More precisely, in the scalar case, some delicate spectral estimates for the seld-adjoint operator ${\mathcal E}''(r)$ were obtained in~\cite{weinsteinMS,weinsteinCpam}, allowing
to get the following energy convexity result.

\begin{theorem}\label{stabilita-scalar}
Le $r$ be a ground state solution of equation~\eqref{scalareqq} with $p<2/N$.
Let $\phi\in H^1(\R^N,\C)$ be such that $\|\phi\|_{2}=\|r\|_{2}$
and define the positive number
\begin{equation*}
\Gamma_\phi=\inf_{\underset{y\in\R^N}{\theta\in[0,2\pi)}}
\|\phi(\cdot)-e^{i\theta}r(\cdot-y))\|_{H^1}^{2}.
\end{equation*}
Then there exist two positive constants $\cala$ and
$C$ such that
\begin{equation*}
\Gamma_\phi\leq C(\E(\phi)-\E(R)),
\end{equation*}
provided that $\E(\phi)-\E(R)<\cala$.
\end{theorem}

For systems, we consider the following definition.

\begin{definition}
	\label{admissible}
	We say that a ground state solution $R=(r_1,r_2)$ of system~\eqref{probel0}
	is admissible for the modulational stability property to hold,
	and we shall write that $R\in {\mathcal R}$, if $r_i\in H^2(\R^N)$ are radial, $|x|r_i\in L^2(\R^N)$, 
    the corresponding solution $\phi_i(t)$ belongs to $H^2(\R^N)$ for all times $t>0$	
	and the following property holds:
	let $\Phi\in \H$ be such that $\|\Phi\|_{2}=\|R\|_{2}$
and define the positive number
\begin{equation}
		\label{defGammaphi}
\Gamma_\Phi:=\inf_{\underset{y\in\R^N}{\theta_1,\theta_2\in[0,2\pi)}}
\|\Phi(\cdot)-(e^{i\theta_1}r_1(\cdot-y),e^{i\theta_2}r_2(\cdot-y))\|_{\H}^{2}.
\end{equation}
Then there exist a continuous function $\rho:\R^+\to\R^+$ with $\frac{\rho(\xi)}{\xi}\to 0$ as $\xi\to 0^+$
and a positive constant $C$ such that
$$
\rho(\Gamma_\Phi)+\Gamma_\Phi\leq C(\E(\Phi)-\E(R)).
$$
In particular, there exist two positive constants $\cala$ and
$C'$ such that
\begin{equation}
		\label{strongconclusion}
\Gamma_\Phi\leq C'(\E(\Phi)-\E(R)),
\end{equation}
provided that $\Gamma_\Phi<\cala$.
\end{definition}

In the one dimensional case, for an important physical class, there exists a ground 
state solution of system~\eqref{probel0} which belongs to the class ${\mathcal R}$ (see~\cite{mps-modstabil}).

\begin{theorem}
Assume that $N=1$, $p\in [1,2)$ and $\beta>1$. Then there exists a ground 
state solution $R=(r_1,r_2)$ of system~\eqref{probel0} which belongs to the class ${\mathcal R}$.
\end{theorem}

\smallskip

\subsection{Statement of the main results}
The action of external forces is represented by a potential $V:\RN\to \R$ satisfying
\beq\label{ipoV}
V\; \text{is a $C^{3}$ function bounded with its derivatives},
\eeq
and we will study the asymptotic behavior (locally in time) as $\eps\to 0$ of the solution of the
following Cauchy problem
\begin{equation}
\label{Part-Sist}
\tag{$S_\eps$}
\begin{cases}
\dys i\vep \partial_{t}\phi_{1}+\frac{\vep^{2}}{2}\Delta \phi_{1}
- V(x)\phi_{1}+\phi_{1}(|\phi_{1}|^{2p}
+\b|\phi_{2}|^{p+1}|\phi_{1}|^{p-1})=0 
& \text{in $\R^N\times\R^+$}, 
\\
\noalign{\vskip4pt}
\dys i\vep \partial_{t}\phi_{2}+\frac{\eps^2}2\Delta \phi_{2}
 - V(x)\phi_{2}+\phi_{2}(|\phi_{2}|^{2p}
+\b|\phi_{1}|^{p+1}|\phi_{2}|^{p-1})=0
 & \text{in $\R^N\times\R^+$},
\\
\noalign{\vskip4pt}
\dys \phi_{1}(x,0)=r_{1}\Big(\frac{x-\tilde x}{\vep}\Bigr)e^{\frac{{\rm i}}{\vep} x\cdot \tilde{\xi}},\qquad
\dys \phi_{2}(x,0)=r_{2}\Bigl(\frac{x-\tilde x}{\vep}\Bigr)e^{\frac{{\rm i}}{\vep} x\cdot \tilde{\xi}},
\end{cases}
\end{equation}
where $\tilde x,\tilde\xi\in\R^N$ $N\geq 1$, the exponent $p$ is such that
\beq\label{ipop}
0<p<2/N
\eeq
It is known (see \cite{fm}) that, under these assumptions, and for any initial
datum in $L^{2}$, there exists a unique  solution $\Phi^{\vep}=(\phi_{1}^{\vep},\phi_{2}^{\vep})$ of the Cauchy problem that exists globally in time.
We have chosen as initial data a scaling of a real vector $R=(r_{1},r_{2})$ belonging to ${\mathcal R}$.
\vskip6pt

The first main result is the following
\bte\label{mainth3}
Let $R=(r_1,r_2)$ be a ground state solution of~\eqref{probel0} which belongs to the class ${\mathcal R}$.
Under assumptions~\eqref{ipoV},~\eqref{ipop}, let  
$\Phi^{\vep}=(\phi_{1}^{\vep},\phi_{2}^{\vep})$ be the family  
of solutions to system \eqref{Part-Sist}.
Furthermore, let $(x(t),\xi(t))$ be the solution of the Hamiltonian system
\begin{equation}
\label{drivingHS}
\begin{cases}
\, \dot{x}(t)=\xi(t)
\\
\, \dot{\xi}(t)=-\nabla V(x(t))
\\
x(0)=\tilde x
\\
\xi(0)=\tilde{\xi}.
\end{cases}
\end{equation}
Then, there exists a locally uniformly bounded family of functions 
$\teta^{\vep}_{i}:\R^+\to S^{1}$, $i=1,2$, such that, 
defining the vector $Q_\eps(t)=(q_{1}^{\vep}(x,t),q_{2}^{\vep}(x,t))$ by
$$
q_{i}^{\vep}(x,t)=r_i\!\!\left(\!\frac{x-x(t)}{\vep}\!\right)
\!e^{\frac{{\rm i}}{\vep}[x\cdot \xi(t)+\teta^{\vep}_{i}(t)]},
$$
it holds 
\begin{equation}
\label{mainconclus}
\|\Phi^{\vep}(t)-Q_\eps(t)\|_{\H_{\vep}}\leq {\mathcal O}(\vep),\quad \text{as $\eps\to0$}
\end{equation}
locally uniformly in time.
\ete

Roughly speaking, the theorem states that, in the semiclassical regime, the modulus of the solution $\Phi^{\vep}$
is approximated, locally uniformly in time, by the admissible real ground state 
$(r_{1},r_{2})$ concentrated in $x(t)$, up to a suitable phase rotation. 
Theorem~\ref{mainth3} can also be read as a description of the {\em slow dynamic} 
of the system close to the invariant manifold of the standing waves generated 
by ground state solutions. This topic has been studied, for the single
equation, in \cite{pw}.

\begin{remark}
Suppose that $\tilde\xi=0$ and $\tilde x$ is a critical point of the potential $V$. Then
the constant function $(x(t),\xi(t))=(\tilde x,0)$, for all $t\in\R^+$, is the solution to system~\eqref{drivingHS}. 
As a consequence, from Theorem~\ref{mainth3}, the approximated solutions is of the form
$$
r_{i}\Big(\!\frac{x-\tilde x}{\vep}\Big)
e^{\frac{{\rm i}}{\vep}\teta^{\vep}_{i}(t)},\quad x\in\R^N,\, t>0,
$$
that is, in the semiclassical regime, the solution concentrates around 
the critical points of the potential $V$.
This is a remark related to ~\cite{mps} where we have considered as
initial data ground states solutions of an associated nonautonomous elliptic problem.
\end{remark}

\begin{remark}
As a corollary of Theorem~\ref{mainth3} we point out that, in the particular case 
of a constant potential, the approximated solution has components 
$$
r_{i}\Big(\frac{x-\tilde x-\tilde\xi t}{\vep}\Big)
e^{\frac{{\rm i}}{\vep}[x\cdot \tilde\xi+\teta^{\vep}_{i}(t)]},
\quad x\in\R^N,\, t>0.
$$
Hence, the mass center $x(t)$ of $\Phi(t,x)$ moves with
constant velocity $\tilde{\xi}$ realizing a uniform motion.
This topic has been tackled, for the single equation, in \cite{weinsteinMS}.
\end{remark}

\begin{remark}
For values of $\beta>1$ both components of the ground states $R$ are nontrivial and,
for $R\in {\mathcal R}$, the solution of the Cauchy problem are approximated by a vector with both nontrivial components. We expect that ground state solutions for $\beta>1$ are unique (up to
translations in $\R^N$) and nondegenerate.
\end{remark}

\vskip2pt
\noindent
We can also analyze the behavior of total momentum density defined by 
\beq\lbl{totmondens}
P^{\vep}(x,t):=p^{\vep}_{1}(x,t)+p^{\vep}_{2}(x,t),\,\quad 
\text{for $x\in\R^N$,\, $t>0$},
\eeq
where
\beq\label{defpeps}
p^{\vep}_{i}(x,t):=\frac1{\vep^{N-1}}
\im\big(\overline{\phi}_{i}^{\vep}(x,t)\nabla\phi_{i}^{\vep}(x,t)\big),\quad
\text{for $i=1,2$,\, $x\in\R^N$,\, $t>0$}.
\eeq
Moreover, let $M(t):=(m_{1}+m_{2})\xi(t)$ be the total momentum of the 
particle $x(t)$ solution of \eqref{drivingHS}, where
\beq\lbl{defmi}
m_{i}:=\|r_{i}\|_{2}^{2},\quad \text{for $i=1,2$}.
\eeq
The information about the asymptotic behavior of $P^{\eps}$ and of the mass densities  $|\phi^{\vep}_{i}|^2/\vep^{N}$ are contained in the following result.

\bte\label{mainth-Bis}
Under the assumptions of Theorem~\ref{mainth3}, there exists $\vep_{0}>0$ such that
\begin{gather*}
\big\|(|\phi^{\vep}_{1}|^{2}/\vep^{N}dx,|\phi^{\vep}_{2}|^{2}/\vep^{N}dx)-\!(m_{1},m_{2})\delta_{x(t)}\big\|_{(C^{2}\times C^{2})^{*}}\leq \calo(\vep^{2}), \\
\noalign{\vskip2pt}
\big\|P^{\vep}(t,x)dx-M(t)\delta_{x(t)}\big\|_{(C^{2})^{*}}
\leq \calo(\vep^{2}),
\end{gather*}
for every $\vep\in(0,\vep_{0})$ and locally uniformly in time. 
\ete

\begin{remark}
Essentially, the theorem states that, in the semiclassical regime, the mass densities of the components $\phi_i$ of the solution $\Phi^{\vep}$ behave 
as a point  particle located in $x(t)$ of mass respectively $m_i$ and the 
total momentum behaves like $M(t)\delta_{x(t)}$. It should be stressed that 
we can obtain the asymptotic behavior for each single mass density, while we can only afford the same result for the total momentum.
The result will follow by a more general
technical statement (Theorem~\ref{mainth2Bis}).
\end{remark}

\begin{remark}
The hypotheses on the potential $V$ can be slightly 
weakened. Indeed, we can assume that  $V$ is bounded from below and that 
$\partial^{\alpha} V$ are bounded only for 
$|\alpha|=2$ or $|\alpha|=3$. This allows to include the important class of 
harmonic potentials (used e.g.\ in Bose-Einstein theory), such as 
$$
V(x)=\frac{1}{2}\sum_{j=1}^N\omega_j^2x_j^2,\quad \omega_j\in\R,\,\, j=1,\dots,N.
$$
Hence, equation~\eqref{drivingHS} reduces to the system of harmonic oscillators
\begin{equation}
\label{systBEc}
\ddot x_j(t)+\omega_j^2x_j(t)=0,\quad j=1,\dots,N.
\end{equation}
For instance, in the 2D case, renaming
$x_1(t)=x(t)$ and $x_2(t)=y(t)$ the ground states solutions are driven 
around (and concentrating) along the lines of a {\em Lissajous curves} 
having {\em periodic} or {\em quasi-periodic} behavior depending on the
case when the ratio $\omega_i/\omega_j$ is, respectively, a {\em rational} 
or an {\em irrational} number. See Figures~\ref{f:esempio5x72x12}
and~\ref{f:esempio5x72x36} below for the corresponding phase portrait 
in some 2D cases, depending on the values of $\omega_i/\omega_j$.
\begin{figure}[h!!!]
\begin{center}
   \includegraphics[scale=.53]{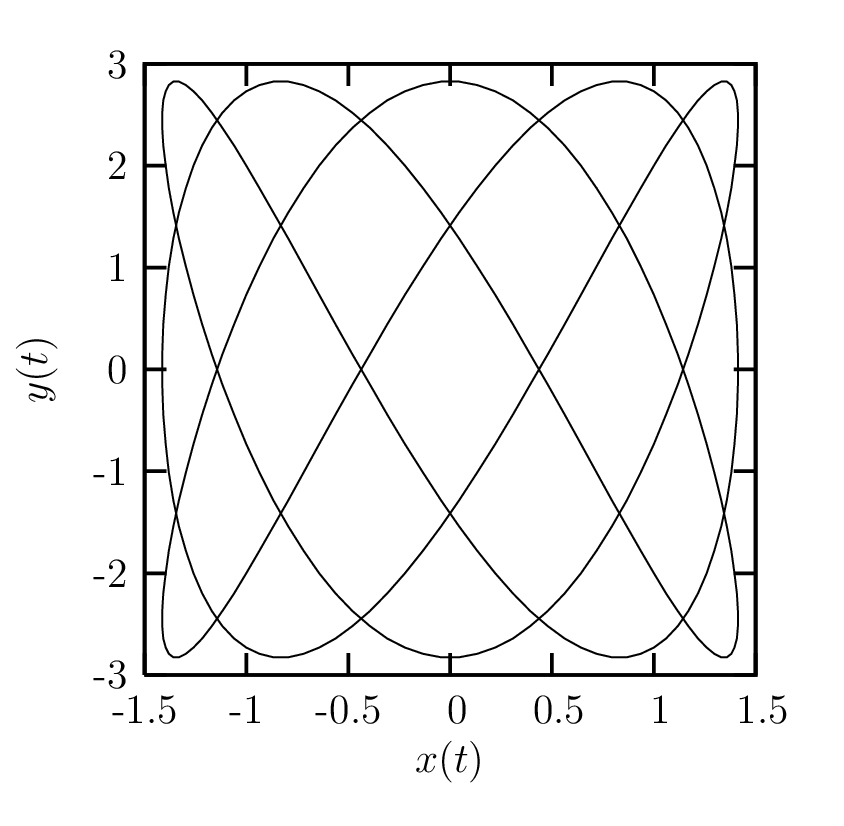}
   \hspace{20pt}
   \includegraphics[scale=.53]{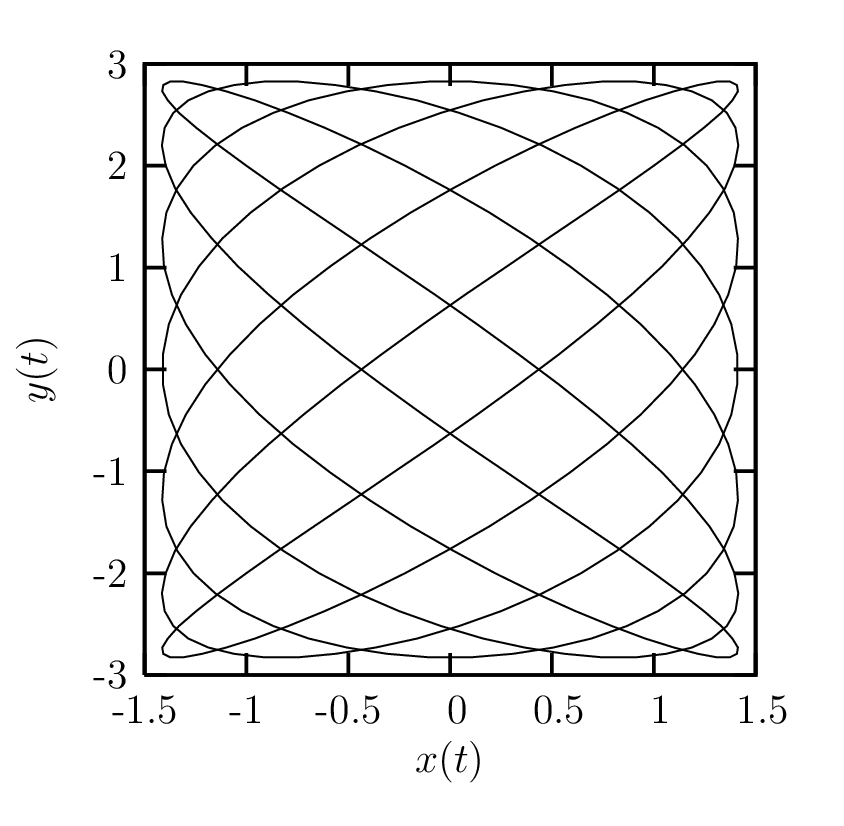}
\end{center}
\caption{Phase portrait of system~\eqref{systBEc} in 2D with $\omega_1/\omega_2=3/5$~(left) and 
	 $\omega_1/\omega_2=7/5$~(right). Notice the periodic behaviour.}
	 \label{f:esempio5x72x12}
\end{figure}
\begin{figure}[h!!!]
\begin{center}
   \includegraphics[scale=.53]{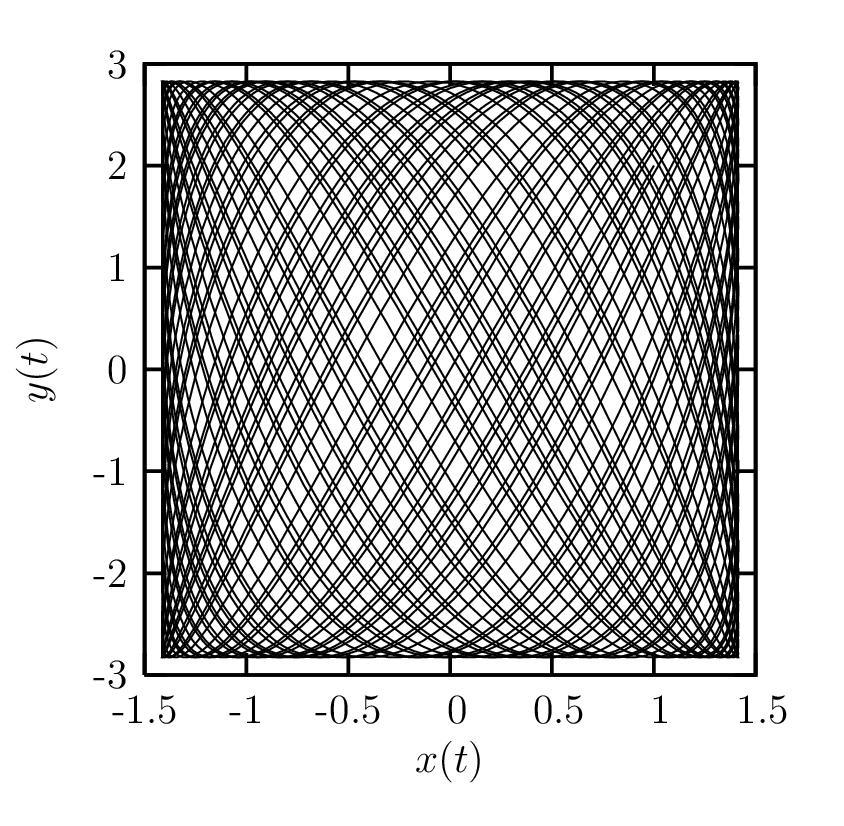}
   \hspace{20pt}
   \includegraphics[scale=.53]{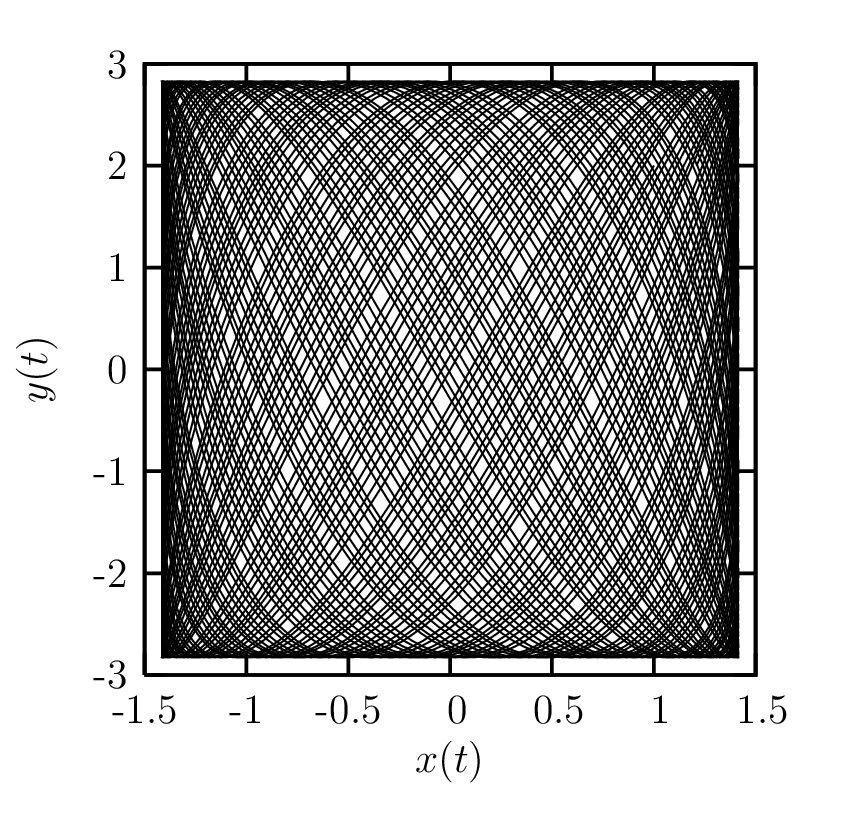}
\end{center}
\caption{Phase portrait of system~\eqref{systBEc} in 
2D with $\omega_1/\omega_2=\sqrt 3/3$
increasing the integration time from $t\in[0,40\pi]$ 
(left) to $t\in[0,60\pi]$ (right). 
Notice the quasi-periodic behaviour, the plane is filling up.}
	 \label{f:esempio5x72x36}
\end{figure}
\end{remark}

\vskip6pt
The paper is organized as follows. 
\vskip2pt
\noindent
In Section~\ref{ingredients} we set up the main ingredients for the proofs as well as state
two technical approximation results (Theorems~\ref{mainth2},~\ref{mainth2Bis}) in a general framework.
In Section~\ref{preli} we will collect some preliminary technical facts that will be useful to prove the results.
In Section~\ref{proofs} we will include the core computations regarding energy and momentum estimates in the
semiclassical regime.
Finally, in Section~\ref{proofMR}, the main results (Theorems~\ref{mainth3} and \ref{mainth-Bis}) will be proved.
\smallskip

\section{A more general Schr\"odinger system}
\label{ingredients}
 
In the following sections we will study the behavior, 
for sufficiently small $\vep$, of a solution 
$\Phi=(\phi_{1},\phi_{2})$ of the more general Schr\"odinger system 
\beq
\label{prob0}
\tag{$F_\eps$}
\begin{cases}
\dys i\vep \partial_{t}\phi_{1}+\frac{\vep^{2}}{2}\Delta \phi_{1}
- V(x)\phi_{1}+\phi_{1}(|\phi_{1}|^{2p}
+\b|\phi_{2}|^{p+1}|\phi_{1}|^{p-1})=0 
& \text{in $\R^N\times\R^+$}, 
\\
\noalign{\vskip4pt}
\dys i\vep \partial_{t}\phi_{2}+\frac{\eps^2}2\Delta \phi_{2}
 - W(x)\phi_{2}+\phi_{2}(|\phi_{2}|^{2p}
+\b|\phi_{1}|^{p+1}|\phi_{2}|^{p-1})=0
 & \text{in $\R^N\times\R^+$},
\\
\noalign{\vskip4pt}
\dys \phi_{1}(0,x)=r_{1}\Big(\frac{x-\tilde x}{\vep}\Big)e^{\frac{{\rm i}}{\vep} x\cdot \tilde{\xi}_1}\qquad
\dys \phi_{2}(0,x)=r_{2}\Big(\frac{x-\tilde x}{\vep}\Big)e^{\frac{{\rm i}}{\vep} x\cdot \tilde{\xi}_2},
\end{cases}
\end{equation}
where $p$ verifies \eqref{ipop}, the potentials $V,\,W$ both satisfy 
\eqref{ipoV} and $(r_{1},r_{2})$ is a real ground state solution of problem \eqref{probel0}. As for the case of a single potential, we get a
unique globally defined $\Phi^{\vep}=
(\phi_{1}^{\vep},\phi_{2}^{\vep})$
that depends continuously on the initial data (see, e.g. \cite[Theorem 1]{fm}). Moreover, if the initial data are chosen in $H^{2}\times H^{2}$, 
then $\Phi^{\vep}(t)$ enjoys the same regularity property for all positive times $t>0$ (see e.g.\ \cite{cazenave}).

\begin{remark}
With no loss of generality, we can assume $V,W\geq 0$. Indeed, if $\phi_1,\phi_2$ is a solution to \eqref{prob0},
since $V,W$ are bounded from below by~\eqref{ipoV}, there exist $\mu>0$ such that $V(x)+\mu\geq 0$ and $W(x)+\mu\geq 0$, for all
$x\in\R^N$. Then $\hat\phi_1=\phi_1e^{-{\rm i}\frac{\mu t}{\vep}}$ and $\hat\phi_2=\phi_2e^{-{\rm i}\frac{\mu t}{\vep}}$ is a solution
of~\eqref{prob0} with $V+\mu$ (resp.\ $W+\mu$) in place of $V$ (resp.\ $W$).
\end{remark}
We will show that the dynamics of $(\phi_{1}^\vep,\phi_{2}^\vep)$
is governed by  the solutions
$$
X=(x_{1},x_{2}):\R\to\R^{2N},\qquad
\Xi=(\xi_{1},\xi_{2}):\R\to\R^{2N},
$$ 
of the following Hamiltonian systems 
\beq
\label{Ham}
\tag{$H$}
\begin{cases}
\dot{x}_{1}(t)=\xi_{1}(t)
\\
\dot{\xi}_{1}(t)=-\nabla V(x_{1}(t))
\\
(x_{1}(0),\xi_{1}(0))=(\tilde x,\tilde{\xi}_{1}),
\end{cases}
\quad
\begin{cases}
\dot{x}_{2}(t)=\xi_{2}(t)
\\ 
\dot{\xi}_{2}(t)=-\nabla W(x_{2}(t))
\\
(x_{2}(0),\xi_{2}(0))=(\tilde x,\tilde{\xi}_{2}).
\end{cases}
\eeq
Notice that the Hamiltonians  related to these systems are
\beq
H_{1}(t)=\frac12 |\xi_{1}(t)|^{2}+V(x_{1}(t)),
\qquad H_{2}(t)=\frac12|\xi _{2}(t)|^{2}+W(x_{2}(t))
\eeq
and are conserved in time. Under assumptions \eqref{ipoV} 
it is immediate to check that the Hamiltonian systems \eqref{Ham} have global solutions.
With respect to the  asymptotic behavior of the solution of \eqref{prob0}
we can prove the following results.

\subsection{Two more general results}

We now state two technical theorems that will yield, as a corollary, Theorems~\ref{mainth3} and~\ref{mainth-Bis}.
\bte\label{mainth2}
Assume \eqref{ipop} and that $V,\,W$ both satisfy \eqref{ipoV}.
Let $\Phi^{\vep}=(\phi_{1}^{\vep},\phi_{2}^{\vep})$ be the family 
of solutions to system~\eqref{prob0}.
Then, there exist $\eps_0>0$, $T^{\eps}_{*}>0$,
a family of continuous functions
$\varrho^\vep:\R^+\to\R$ with $\varrho^\vep(0)=\calo(\vep^{2})$, locally uniformly bounded sequences 
of functions $\teta^{\vep}_{i}:\R^+\to S^{1}$ and a positive constant
$C$, such that, defining the vector 
$Q_\vep(t)=(q_1^\vep(x,t),\,q_2^\vep(x,t))$ by 
$$
q_i^\vep(x,t)=r_i\!\!\left(\!\frac{x-x_{1}(t)}{\vep}\!\right)
\!e^{\frac{{\rm i}}{\vep}[x\cdot \xi_{i}(t)+\teta^{\vep}_{i}(t)]},\quad i=1,2
$$
it results 
$$
\|\Phi^{\vep}(t)-Q_\vep(t)\|_{\H_{\vep}}\leq
C\sqrt{\varrho^\vep(t)+\left(\frac{\varrho^\vep(t)}{\eps}\right)^2},
$$
for all $\eps\in (0,\eps_0)$ and all $t\in [0,T^{\eps}_{*}]$,
where $x_1(t)$ is the first component of the Hamiltonian system for $V$ in~\eqref{Ham}.
\ete

\begin{remark}
Theorem~\ref{mainth2} is quite instrumental in the context of our paper, 
as we cannot guarantee in the general case of different potentials that the function $\varrho^\vep$ is small as $\vep$ vanishes, locally uniformly
in time. Moreover, the time dependent shifting of the components $q_i$ into $x_1(t)$ is quite arbitrary, a similar statement could
be written with the component $x_2(t)$ in place of $x_1(t)$, this arbitrariness
is a consequence of the same initial data $\tilde x$ in~\eqref{Ham} for
both $x_{1}$ and $x_{2}$. The task of different initial data in~\eqref{Ham}
for $x_{1}$ and $x_{2}$ is to our knowledge an open problem.
\end{remark}

In the following, if $\xi_i$ are the second components of the systems in \eqref{Ham}, we set
\beq\lbl{totmonPart}
M(t):=m_{1}\xi_{1}(t)+m_{2}\xi_{2}(t),\,\quad t>0.
\eeq
If $\Phi^{\vep}=(\phi_{1}^{\vep},\phi_{2}^{\vep})$ is the family of 
solutions to~\eqref{prob0}, we have the following

\bte\label{mainth2Bis}
There exist $\vep_{0}>0$ and $T^{\vep}_{*}>0$ and a family of continuous functions
$\varrho^\vep:\R^+\to\R$ with $\varrho^\vep(0)=\calo(\vep^{2})$ such that
$$
\begin{array}{cc}
\big\|(|\phi^{\vep}_{1}|^{2}/\vep^{N}dx,|\phi^{\vep}_{2}|^{2}/\vep^{N}dx)-\!(m_{1},m_{2})\delta_{x_{1}(t)}\big\|_{(C^{2}\times C^{2})^{*}}\leq \varrho^\vep(t),
\\\\
\big\|P^{\vep}(t,x)dx-M(t)\delta_{x_{1}(t)})\big\|_{(C^{2})^{*}}
\leq \varrho^\vep(t),
\end{array}
$$
for every $\vep\in(0,\vep_{0})$ and all $t\in[0,T^{\vep}_{*}]$.
\ete
\section{Some preliminary results} \lbl{preli}

In this section we recall and show some results we will use in proving
Theorems \ref{mainth3}, \ref{mainth-Bis}, \ref{mainth2} and \ref{mainth2Bis}.
First we recall the following conservation laws.
\begin{proposition}\lbl{conse}
The mass components of a solution $\Phi$ of \eqref{prob0},
\beq\lbl{masse}
{\mathcal N}^{\vep}_{i}(t):=\frac1{\vep^{N}}\|\phi^{\vep}_{i}(t)\|_{\elle2}^{2},\quad \text{for $i=1,2$,\,\, $t>0$},
\eeq
are conserved in time. Moreover, also the total energy defined by
\beq\label{consE}
 E^{\vep}(t)=E_{1}^{\vep}(t)+E_{2}^{\vep}(t)
\eeq
is conserved as time varies, where
\begin{align*}
E_{1}^{\vep}(t) &=\frac1{2\vep^{N-2}}
\|\nabla\phi^{\vep}_{1}\|_{\elle2}^{2}
+\frac1{\vep^{N}}\int  V(x)|\phi^{\vep}_{1}|^{2}dx 
-\frac1{2\vep^{N}}\int F_{\beta}(\Phi^{\vep})dx, \\
E_{2}^{\vep}(t)&=\frac1{2\vep^{N-2}}\|\nabla\phi^{\vep}_{2}\|_{2}^{2}
+\frac1{\vep^{N}}\int W(x)|\phi^{\vep}_{2}|^{2}dx
-\frac1{2\vep^{N}}\int F_{\beta}(\Phi^{\vep})dx.
\end{align*}
\end{proposition}
\bdim
This is a standard fact. For the proof, see e.g.\ \cite{fm}.
\edim
\vskip6pt
\bos
From the preceding proposition we obtain that, due to the form of our initial data,
 the mass components ${\mathcal N}_{i}^{\vep}(t)$ do not actually depend on $\vep$. Indeed, for $i=1,2$,
\beq\label{eqmi}
{\mathcal N^{\vep}_{i}(t)}={\mathcal N^{\vep}_{i}(0)}=\frac1{\vep^{N}}\int|\phi^{\vep}_{i}(x,0)|^{2}dx=\frac1{\vep^{N}}\int \Big|r_{i}\Big(\frac{x-\tilde x}{\vep}\Big)
\Big|^{2}dx=m_{i}.
\eeq
Thus, the quantities $\phi^{\vep}_{i}/\vep^{N/2}$ have constant norm in $L^2$ equal, respectively, 
to $m_{i}$. In Theorem~\ref{mainth2Bis} we will show that, for sufficiently small values of $\vep$, 
the mass densities behave, point-wise with respect to $t$,
as a $\delta$ functional concentrated in $x_{1}(t)$.
\eos
In the following we will often make use of the following simple Lemma.
\begin{lemma}\label{pote}
Let $A\in C^{2}(\rn)$ be such that $A,D_jA, D^2_{ij}A$ are uniformly 
bounded and let $R=(r_{1},r_{2})$ 
be a ground state solution of problem~\ref{probel0}. Then,
for every $y\in \rn$ fixed, there exists a positive constant
$C_{0}$ such that
\beq
\left|\int \left[A(\vep x+y)-A(y)\right] r_{i}^{2}(x)dx\right|\leq C_{0}\vep^{2}.
\eeq
\end{lemma}
\bdim
By virtue of the regularity properties of the function $A$ and Taylor expansion Theorem we get
\begin{align*}
\frac1{\vep^{2}}\left|\int\left[A(\vep x+y)-A(y)\right]r_{i}^{2}(x)dx\right|
&\leq
\frac1{\vep}|\nabla A(y)|\left|\int x r_{i}^{2}(x)dx\right|
\\
&+\|{\rm Hes}(A)\|_{\infty}\int |x|^{2}r_{i}^{2}(x)dx
\end{align*}
where $\|{\rm Hes}(A)\|_{\infty}$ denotes the $L^{\infty}$ norm of the Hessian
matrix associated to the function $A$. The first integral on the right 
hand side is zero since each component $r_i$ is radial. The
second integral is finite, since $|x|r_{i}\in L^2(\R^N)$.
\edim
\vskip6pt
In order to show the desired asymptotic behavior 
we will use the following property of the functional $\delta_{y}$ on the space $C^{2}(\rn)$.
\begin{lemma}\lbl{propdelta}
There exist $K_{0},\,K_{1},\,K_{2}$ positive constants, such that, if 
$\|\delta_{y}-\delta_{z}\|_{C^{2*}}$$\leq K_{0}$ then
$$
K_{1}|y-z|\leq \|\delta_{y}-\delta_{z}\|_{C^{2*}}\leq K_{2}|y-z|
$$
\end{lemma}
\bdim
For the proof see \cite[Lemma 3.1, 3.2]{Keerani2}.
\edim
\vskip6pt
The following lemma will be used in proving our main result.
\begin{lemma}\label{propalfa}
Let $\Phi^\vep=(\phi_{1}^\vep,\phi_{2}^\vep)$ be a solution of \eqref{prob0} and
consider the vector functions $\alpha_{i}:\R\to \rn$ defined by
\beq
\label{alfai}
\alpha^{\vep}_{i}(t)=\int p^{\vep}_{i}(x,t)dx-m_{i}\xi_{i}(t),\quad t>0,\, i=1,2,
\eeq
where the $\xi_{i}$s are defined in~\eqref{Ham} and the $m_{i}$s
are defined in~\eqref{defmi}, for $i=1,2$. Then $\{t\mapsto \alpha^{\vep}_{i}(t)\}$ is 
a continuous function and $\alpha^{\vep}_{i}(0)=0$, for $i=1,2$. 
\end{lemma}
\bos
The integral in \eqref{alfai} defines a vector whose components 
are the integral of $\im(\overline{\phi^{\vep}_{i}}\partial\phi^{\vep}_{i}/\partial x_{j})/\vep^{N-1}$ for $j=1,\dots,N$, so that $\alpha_{i}^{\vep}:\R\to\rn.$
\eos
\bdim
The continuity of $\alpha_{i}$ immediately follows from the
regularity properties of the solution $\phi^{\vep}_{i}$.
In order to complete the proof, first note that, for all $x\in\R^N$,
$$
\bar \phi_{i}^{\vep}(x,0)\nabla \phi_{i}^{\vep}(x,0)=
\frac{{\rm i}}{\vep}\tilde\xi_{i}r_{i}^{2}\Big(\frac{x-\tilde x}\vep\Big)
+\frac1{\vep}r_{i}\Big(\frac{x-\tilde x}\vep\Big)\nabla r_{i}
\Big(\frac{x-\tilde x}\vep\Big),
$$
so that, as $r_{i}$ is a real function, the conclusion follows by a change of variable.
\edim
\vskip6pt
\begin{lemma}\label{propeta}
Let $V$ and $W$ both satisfying assumptions~\eqref{ipoV} and let
$\Phi^\vep=(\phi_{1}^\vep,\phi_{2}^\vep)$ be  a solution of~\eqref{prob0}. Moreover,
let $A$ a positive constant defined by
\beq\label{defa}
A=K_{1}\sup_{[0,T_{0}]}\left[|x_{1}(t)|+|x_{2}(t)|\right]+K_{0}
\eeq
where $x_{i}(t)$ is defined in \eqref{Ham}, $K_{0}$ and $K_1$
are defined in Lemma~\ref{propdelta}, and  let $\chi$ be a 
$C^{\infty}(\rn)$ function such that $0\leq\chi\leq 1$ and 
\beq\lbl{defchi}
\chi(x)=1\qquad\text{if $|x|< A$},\qquad\quad \chi(x)=0\qquad
\text{if $|x|> 2A$}. 
\eeq
Then the  functions
\beq\label{eta}
\left\{
\begin{array}{ll}
\dys \eta^{\vep}_{1}(t)=\dys m_{1}V(x_{1}(t))-
\frac1{\vep^{N}}\int \chi(x)V(x)|\phi^{\vep}_{1}(x,t)|^{2}dx,&
\medskip\\
\dys \eta^{\vep}_{2}(t)=m_{2}W(x_{2}(t))-
\frac1{\vep^{N}}\int \chi(x)W(x)|\phi^{\vep}_{2}(x,t)|^{2}dx.&
\end{array}
\right.
\eeq
are continuous  and satisfy
 $\left|\eta^{\vep}_{i}(0)\right|=\calo(\vep^{2})$ for $i=1,2$. 
\end{lemma}
\bdim
The continuity of $\eta_{i}^\vep$ immediately follows from the
regularity properties of the solution $\phi^{\vep}_{i}$.
We will prove the conclusion only for $\eta^{\vep}_{1}(0)$, the result for 
$\eta^{\vep}_{2}(0)$ can be showed in an analogous way. We have
\begin{align*}
 \left|\eta^{\vep}_{1}(0)\right|&=\left|m_{1}V(x_{1}(0))-\frac1{\vep^{N}}\int \chi(x)V(x)|\phi^{\vep}_{1}(x,0)|^{2}dx\right|
\\
&\leq \left|m_{1}V(\tilde x)-\frac1{\vep^{N}}\int V(x)r_{1}^{2}\Big(\frac{x-\tilde x}\vep\Big)dx\right|
\\
&+\frac1{\vep^{N}}\int_{|x|>A}\!\!\left(1-\chi(x)\right)V(x)
r_{1}^{2}\Big(\frac{x-\tilde x}\vep\Big)dx.
\end{align*}
Then, by Lemma~\ref{pote}, and a change of variables imply
\begin{align*}
\left|\eta^{\vep}_{1}(0)\right|&\leq \calo(\vep^{2})+
\int\!\!\left(1-\chi(\tilde x+\vep y)\right)V(\tilde x
+\vep y)r_{1}^{2}\left(y\right)dy.
\end{align*}
The properties of $\chi$ and $r_1$ and assumption \eqref{ipoV}  yield the conclusion.
\edim
\vskip6pt
We will also use the following identities.
\begin{lemma}\lbl{ideder}
The following identities holds for $i=1,2$.
\beq\lbl{derfi}
\frac1{\vep^{N}}\frac{\partial |\phi^{\vep}_{i}|^{2}}{\partial t}(x,t)=-{\rm div}_{x}\, p^{\vep}_{i}(x,t),\quad
x\in\R^N,\, t>0.
\eeq
Moreover, for all $t>0$, it results
\begin{equation}
\lbl{derpot}
\int \frac{\partial P^{\vep}}{\partial t}(x,t)dx=-\frac1{\vep^{N}}\int\nabla V(x)|\phi^{\vep}_{1}(x,t)|^{2}dx
-\frac1{\vep^{N}}\int\nabla W(x)|\phi^{\vep}_{2}(x,t)|^{2}dx,
\end{equation}
where $P^\vep(x,t)$ is the total momentum density defined in \eqref{totmondens}.
\end{lemma}

\begin{remark}
It follows from identity~\eqref{derpot} that for systems with constant potentials the total momentum $\int P^{\vep}dx$
is a constant of motion.
\end{remark}

\begin{remark}
As evident from identity~\eqref{derpot} as well as physically reasonable, in the case of systems of Schr\"odinger equations, the balance for the momentum 
needs to be stated for the sum $P^{\vep}$ instead on the single components $p^{\vep}_{i}$. See also 
identities~\eqref{p1control} and~\eqref{p2control} in the proof, where the coupling terms appear.
\end{remark}

\bdim
In order to prove identity \eqref{derfi} note that
$$
-{\rm div}_{x}p^{\vep}_{i}=-\frac1{\vep^{N-1}}\im(\bar\phi^{\vep}_{i}\Delta\phi_{i}^{\vep}), \qquad 
\frac1{\vep^{N}}\frac{\partial |\phi^{\vep}_{i}|^{2}}{\partial t}
=\frac2{\vep^{N}}\re((\phi^{\vep}_{i})_t\bar\phi^{\vep}_{i})
$$
Since $\phi^{\vep}_{i}$ solves the corresponding equation in
system~\eqref{prob0}, we can multiply the equation by 
$\bar \phi^{\vep}_{i}$ and add this identity to its conjugate; 
the conclusion follows from the properties of the nonlinearity.
Concerning identity \eqref{derpot}, observe first that, setting 
$(p^{\vep}_1)_j(x,t)=\vep^{1-N}\im(\overline{\phi}_1^{\vep}(x,t)\partial_j\phi_1^{\vep}(x,t))$
for any $j$ and $\partial_j=\partial_{x_j}$, it holds
\begin{align*}
\frac{\partial (p^{\vep}_1)_j}{\partial t} & =\vep^{1-N}\im(\partial_t\overline{\phi}_1^{\vep}\partial_j\phi_1^{\vep})
+\vep^{1-N}\im(\overline{\phi}_1^{\vep}\partial_j(\partial_t\phi_1^{\vep})) \\
& =\vep^{1-N}\im(\partial_t\overline{\phi}_1^{\vep}\partial_j\phi_1^{\vep})
+\vep^{1-N}\im (\partial_j\big(\overline{\phi}_1^{\vep}\partial_t\phi_1^{\vep}\big))
- \vep^{1-N}\im (\partial_j\overline{\phi}_1^{\vep}\partial_t\phi_1^{\vep}) \\
\noalign{\vskip3pt}
& =2\vep^{1-N}\im(\partial_t\overline{\phi}_1^{\vep}\partial_j\phi_1^{\vep})
+\vep^{1-N}\im (\partial_j\big(\overline{\phi}_1^{\vep}\partial_t\phi_1^{\vep}\big)).
\end{align*}
In particular the second term integrates to zero. Concerning the first addendum, take the first
equation of system \eqref{prob0}, conjugate it and multiply it by $2\eps^{-N}\partial_j\phi_1$. It follows
\begin{align*}
2\vep^{1-N} \im(\partial_{t}\overline{\phi_{1}}^{\vep}\partial_j\phi_1^{\vep}) &=-\vep^{2-N}\re(\Delta 
\overline{\phi_{1}}^{\vep}\partial_j\phi_1^{\vep})+ 2\eps^{-N}V(x)\re(\overline{\phi_{1}}^{\vep}\partial_j\phi_1^{\vep}) \\
\noalign{\vskip3pt}
& -2\eps^{-N}|\phi_{1}^{\vep}|^{2p}\re(\overline{\phi_{1}}^{\vep}\partial_j\phi_1^{\vep}) -2\b\eps^{-N}|
\phi_{2}^{\vep}|^{p+1}|\phi_{1}^{\vep}|^{p-1}\re(\overline{\phi_{1}}^{\vep}\partial_j\phi_1^{\vep}) \\
\noalign{\vskip3pt}
&=-\vep^{2-N}\re(\partial_i\big(\partial_i \overline{\phi_{1}}^{\vep}\partial_j\phi_1^{\vep}))+\vep^{2-N}\partial_j\Big(\frac{|\partial_i \phi_{1}^{\vep}|^{2}}{2}\Big) \\
\noalign{\vskip3pt}
& + \eps^{-N}\partial_j\left(V(x) |\phi_{1}^{\vep}|^2\right)-\eps^{-N}\partial_jV(x)|\phi_1^\vep|^2  \\
& -\eps^{-N}\partial_j\Big(\frac{|\phi_{1}^{\vep}|^{2p+2}}{p+1}\Big)-2\b\eps^{-N}|\phi_{2}^{\vep}|^{p+1}  \partial_j\Big(\frac{|\phi_{1}^{\vep}|^{p+1}}{p+1}\Big).
\end{align*}
Of course, one can argue in a similar fashion for the second component $\phi_2$.
Then, taking into account that all the terms in the previous identity but 
$\partial_jV(x)|\phi_1^\vep|^2$ and $|\phi_{2}^{\vep}|^{p+1}  \partial_j|\phi_{1}^{\vep}|^{p+1}$ integrate to zero
due to the $H^2$ regularity of $\phi_1$, we reach
\begin{align}
\label{p1control}
\int \frac{\partial (p^{\vep}_1)_j}{\partial t}dx&=-\frac1{\vep^{N}}\int \frac{\partial V}{\partial x_j}(x)|\phi^{\vep}_{1}|^{2}dx
-\frac{2\beta}{\eps^{N}}\int |\phi^{\vep}_{2}|^{p+1} \partial_j\Big(\frac{|\phi_{1}^\vep|^{p+1}}{p+1}\Big)dx        \\
\label{p2control}
\int \frac{\partial (p^{\vep}_2)_j}{\partial t}dx&=-\frac1{\vep^{N}}\int \frac{\partial W}{\partial x_j}(x)   |\phi^{\vep}_{2}|^{2}dx
-\frac{2\beta}{\eps^{N}}\int |\phi^{\vep}_{1}|^{p+1}\partial_j\Big(\frac{|\phi_{2}^\vep|^{p+1}}{p+1}\Big)dx.
\end{align}
Adding these identities for any $j$ and taking into account that by
the  regularity pro\-per\-ties of 
$\phi^{\vep}_{i}$ it holds $\int \partial_j(|\phi_{1}^\vep|^{p+1}|\phi_{2}^\vep|^{p+1})dx=0$,
formula~\eqref{derpot} immediately follows.
\edim

\section{Energy, mass and momentum estimates}\lbl{proofs}

\subsection{Energy estimates in the semiclassical regime}
In order to obtain the desired asymptotic behavior stated in Theorems \ref{mainth3}, \ref{mainth-Bis}, \ref{mainth2} and \ref{mainth2Bis},
we will first prove a key inequality concerning the functional $\E$ defined in \eqref{Ecorsivo}.
As pointed out in the introduction, the main ingredients involved are the conservations laws of the 
Schr\"odinger system and of the Hamiltonians functions and a modulational stability 
property for admissible ground states.

The idea is to evaluate the functional $\E$ on the vector $\Upsilon^{\vep}=(v^{\vep}_{1},v^{\vep}_{2})$ whose components are given by
\beq\label{veps}
v^{\vep}_{i}(x,t) =e^{-\frac{{\rm i}}{\vep}\xi_{i}(t)\cdot [\vep x+x_{1}(t)]}
\,\phi^{\vep}_{i}(\vep x+x_{1}(t),t)
\eeq
where $X=(x_{1},x_{2}),\,\Xi=(\xi_{1},\xi_{2})$ are the solution of the system \eqref{Ham}. More precisely, we will prove the following result.
\begin{theorem}\label{chiave}
Let $\Phi^{\vep}=(\phi_{1}^{\vep},\phi_{2}^{\vep})$ be a family 
of solutions of $(F_{\vep})$, and let $\Upsilon^{\vep}$ be the vector defined in
\eqref{veps}.
Then, there exist $\eps_{0}$ and $T^{\eps}_{*}$ such that 
for every $\eps\in(0,\eps_{0})$ and for every $t\in [0,T^{\eps}_{*})$, 
it holds
\be\label{stima}
0\leq \E(\Upsilon^{\vep})-\E(R)\leq  \alpha^{\vep}+\eta^{\vep}+\calo(\vep^{2}),
\ee
where we have set
\begin{equation}
\label{defnalphaeta}
\alpha^{\vep}(t)=\big|(\xi_{1}(t),\xi_{2}(t))\cdot(\alpha_{1}^{\vep}(t), \alpha_{2}^{\vep}(t))\big|,\quad 
\eta^\vep(t)=|\eta_{1}^\vep(t)+\eta_{2}^\vep(t)|,
\end{equation} 
$\alpha_{i},\,\eta_{i}$ are given in~\eqref{alfai},~\eqref{eta} and $R=(r_{1},r_{2})$ 
is the real ground state belonging to the class ${\mathcal R}$ taken as initial datum in~\eqref{prob0}. 
Moreover, there exist families of functions 
$\theta^{\vep}_{i}$, $y_{1}^{\vep}$ and a positive constant $L$ 
such that
\begin{equation}\lbl{quasifine}
{\Big \|}\Phi^{\vep} -\Big(e^{\frac{{\rm i}}{\vep}(x\xi_{1}+\theta^{\vep}_{1})}
r_{1}\Big(\frac{x-y^{\vep}_{1}}\vep\Big),\,\,
e^{\frac{{\rm i}}{\vep}(x\xi_{2}+\theta^{\vep}_{2})}r_{2}\Big(\frac{x-y^{\vep}_{1}}
\vep\Big)\Big)
{\Big \|}^{2}_{\H_{\vep}}\!\leq  L\big[\alpha^{\vep}+\eta^{\vep}+\calo(\vep^{2})\big],
\end{equation}
for every $\eps\in(0,\eps_{0})$ and all $t\in [0,T^{\eps}_{*})$.
\end{theorem}
\bdim
By a change of variable and Proposition~\ref{conse}, we get 
\begin{equation}\label{normel2}
\dys \|v^{\vep}_{i}(\cdot,t)\|_{2}^{2}=\dys
\|\phi^{\vep}_{i}(\vep x+x_{1}(t),t)\|_{2}^{2}=
\frac1{\vep^{N}}\|\phi^{\vep}_{i}(\cdot,t)\|_{2}^{2}= m_{i},\quad t>0,\, i=1,2,
\end{equation}
where $m_i$ are defined in \eqref{defmi}. Hence the mass of $v^{\vep}_{i}$ is conserved during the evolution.
Moreover, by a change of variable, and recalling definition \eqref{defpeps}
we have
\begin{align*}
\E(\Upsilon^{\vep}) &=
\frac1{2\vep^{N-2}}\|\nabla \Phi^{\vep}\|_{2}^{2} +
\frac12\left(m_{1}|\xi_{1}|^2+m_{2}|\xi_{2}|^2\right)
-\frac1{\vep^{N}}F_{\beta}(\Phi^{\vep}) \\
&-\int (\xi_{1}(t),\xi_{2}(t))\cdot (p_{1}^{\vep}(x,t), p_{2}^{\vep}(x,t))dx.
\end{align*}
Then, taking into account the form of the total energy functional, we obtain
\begin{align*}
\E(\Upsilon^{\vep}) &= E^{\vep}(t) -\frac1{\vep^{N}}
\int\left[ V(x)|\phi^{\vep}_{1}|^{2}+ W(x)|\phi^{\vep}_{2}|^{2}\right]dx
+\frac12\left(m_{1}|\xi_{1}|^2+m_{2}|\xi_{2}|^2\right)
\\
&-\int (\xi_{1}(t),\xi_{2}(t))\cdot (p_{1}^{\vep}(x,t), p_{2}^{\vep}(x,t))dx.
\end{align*}
Moreover, using Proposition~\ref{conse} and performing
a change of variable we get
\begin{align*}
E^{\vep}(t)&=E^{\vep}(0)=E^{\vep}\Big(r_{1}\Big(\frac{x-\tilde{x}}{\vep}\Big)e^{\frac{{\rm i}}{\vep}x\cdot\tilde{\xi}_{1}},r_{2}\Big(\frac{x-
\tilde{x}}{\vep}\Big)e^{\frac{{\rm i}}{\vep}x\cdot\tilde{\xi}_{2}}\Big)
\\
&=\E(R)+\frac12\big(m_{1}|\tilde\xi_{1}|^2+m_{2}|\tilde\xi_{2}|^2\big) \\
&+\int \left[V(\vep x+\tilde x)|r_{1}|^{2}+W(\vep x+\tilde x)|r_{2}|^{2}\right]dx,
\end{align*}
this joint with Lemma~\ref{pote} and the conservation of the Hamiltonians
$H_{i}(t)$ yield
\begin{align*}
\E(\Upsilon^{\vep})-\E(R)=&\frac12\Big[m_{1}(|\tilde\xi_{1}(t)|^2+|\xi_{1}(t)|^2)
+m_{2}(|\tilde\xi_{2}(t)|^2+|\xi_{2}(t)|^2)\Big]  \\
& -\int (\xi_{1}(t),\xi_{2}(t))\cdot (p_{1}^{\vep}(x,t), p_{2}^{\vep}(x,t))dx
\\
\noalign{\vskip3pt}
&\!\!+m_{1}V(\tilde x)+m_{2}W(\tilde x)
-\frac1{\vep^{N}}
\int\big[ V(x)|\phi^{\vep}_{1}|^{2}+ W(x)|\phi^{\vep}_{2}|^{2}\big]dx
\\
\noalign{\vskip3pt}
=& \,m_{1}\big[|\xi_{1}(t)|^{2}+V(x_{1}(t))\big]+
m_{2}\big[|\xi_{2}(t)|^{2}+W(x_{2}(t))\big]
\\
\noalign{\vskip3pt}
&\!\!- \int (\xi_{1}(t),\xi_{2}(t))\cdot (p_{1}^{\vep}(x,t), p_{2}^{\vep}(x,t))dx
\\
\noalign{\vskip3pt}
&\!\!-\frac1{\vep^{N}}
\int\big[ V(x)|\phi^{\vep}_{1}|^{2}+ W(x)|\phi^{\vep}_{2}|^{2}\big]dx
+{\mathcal O}(\vep^{2})
\end{align*}
Using the definitions of $\alpha_{i}$ and $\eta_{i}$, we get
\begin{align*}
\E(\Upsilon^{\vep})-\E(R)&\leq -(\xi_{1}(t),\xi_{2}(t))\cdot
(\alpha_{1}^{\vep}(t), \alpha_{2}^{\vep}(t))+\eta^{\vep}(t)
\\
&\;\;\;-\frac1{\vep^{N}}\int (1-\chi(x))\big[
V(x)|\phi^{\vep}_{1}|^{2}+W(x)|\phi^{\vep}_{2}|^{2}\big]dx+
\calo(\vep^{2}),
\end{align*}
Since  $V$ and $W$ are nonnegative functions, by~\eqref{defnalphaeta} it follows that
$$
\E(\Upsilon^{\vep}(t))-\E(R)\leq \alpha^{\vep}(t)+\eta^{\vep}(t)+\calo(\vep^{2}).
$$
Finally, \eqref{caravar} and \eqref{normel2} imply the first conclusion of
Theorem \ref{chiave}, where the positive time $T^{\vep}_{*}$ is built up as follows.
Let $T_0>0$ (to be chosen later). In order to conclude the proof of the result,
notice that $\alpha_{i}(t)$ and $\eta_{i}(t)$ 
are continuous functions by Lemmas~\ref{propalfa} and~\ref{propeta}. Moreover, let
$\Gamma_{\Upsilon^\eps(t)}$ be the positive number given in~\eqref{defGammaphi}
for $\Phi=\Upsilon^\eps(t)$. Notice that $\{t\mapsto \Gamma_{\Upsilon^\eps(t)}\}$
is continuous and, in view of the choice of the initial data~\eqref{dati}, it holds $\Gamma_{\Upsilon^\eps(0)}=0$.
Hence, for every fixed $h_0,h_1>0$,  we can define the time $T_{*}^{\vep}>0$ by
\beq\label{defteps}
T^{\vep}_{*}:=\sup\{t\in[0,T^{0}]:\,\max\{\alpha^\vep(s),\,
\eta^\vep(s)\}\leq h_{0},\,\,\Gamma_{\Upsilon^\eps(s)}\leq h_1,\,\,\text{for all $s\in(0,t)$}\},
\eeq
Notice that, by~\eqref{stima} and choosing $\vep_{0}$ sufficiently small 
we derive, for all $t\in [0,T^{\vep}_*)$ and $\eps\in(0,\eps_0)$, that
$0\leq \E(\Upsilon^{\vep}(t))-\E(R)\leq 3h_{0}$. Now we choose $h_1$ so small that $h_1< {\mathcal A}$,
where ${\mathcal A}$ is the constant appearing in the statement of the admissible 
ground state (see Definition~\eqref{admissible}). Therefore, from 
conclusion~\eqref{strongconclusion}, there exists a positive constant $L$ such that
\begin{equation}
	\label{Biggammaest}
\Gamma_{\Upsilon^\eps(t)}\leq L(\E(\Upsilon^{\vep}(t))-\E(R))\leq L\big[\alpha^{\vep}(t)+\eta^{\vep}(t)+\calo(\vep^{2})\big],
\end{equation}
for every $t\in[0,T^{\vep}_{*})$ and all $\vep\in (0,\vep_{0})$. 
In turn, there exist two
families of functions $\tilde y^{\vep}(t)$ and $\tilde\theta^{\vep}_{i}(t)$ $i=1,2$ such that
\beq
\|\Upsilon^{\vep}(\cdot,t)-
\big(e^{{\rm i}\tilde \theta^{\vep}_{1}(t)}r_{1}(\cdot+\tilde y^{\vep}(t)),
e^{{\rm i}\tilde \theta^{\vep}_{2}(t)}r_{2}(\cdot+\tilde y^{\vep}(t))
\big)
\|_{\H}^{2}\leq 
L\big[\alpha^{\vep}(t)+\eta^{\vep}(t)+\calo(\vep^{2})\big]
\eeq
for every $t\in[0,T^{\vep}_{*})$ and all $\vep\in (0,\vep_{0})$. 
Making a change of variable and using the notation
$$
\theta^{\vep}_{i}(t):=\vep\tilde \theta^{\vep}_{i}(t),\qquad 
y^{\vep}_{1}(t):=x_{1}(t)-\vep\tilde y^{\vep}(t),
$$
the assertion follows.
\edim
\vskip6pt
\begin{remark}
The previous result holds for every $t\in [0,T^{\eps}_{*})$ where $T^{\eps}_{*}$ is found in \eqref{defteps} and $T^{\eps}_{*}\leq T_{0}$. But,
we have not fixed $T_{0}$ yet. This will be done in Lemma \ref{y}.
\end{remark}

\subsection{Mass and total momentum estimates}
The next lemmas will be used to prove the desired asymptotic behavior.
We start with the study of the asymptotic behavior of the mass densities
and the total momentum density. From now on we shall set
$$
\hat\alpha^\vep(t):=\alpha^\vep(t)+|\alpha^\vep_1(t)+\alpha^\vep_2(t)|,\quad t>0.
$$
\begin{lemma}\label{tecnico1}
There exists a positive constant $L_{1}$ such that 
$$
\begin{array}{cc}
\|(|\phi^{\vep}_{1}|^{2}/\vep^{N}dx,|\phi^{\vep}_{2}|^{2}/\vep^{N}dx)-\!(m_{1},m_{2})\delta_{y^{\vep}_{1}(t)}\|_{(C^{2}\times C^{2})^{*}}
\\\\
+
\|P^{\vep}(t,x)dx-M(t)\delta_{y^{\vep}_{1}(t)}\|_{(C^{2})^{*}}
\\\\
\leq L_{1}\left[\hat\alpha^{\vep}(t)+\eta^{\vep}(t)+\calo(\vep^{2})\right],
\end{array}
$$
for every $t\in[0,T^{\vep}_{*}]$ and $\vep\in(0,\vep_{0})$.
\end{lemma}
\bos
This result will immediately imply Theorem \ref{mainth3}, once we have 
shown the desired asymptotic behavior of $\hat\alpha^{\vep}(t)+\eta^{\vep}(t)$ and of the functional $\delta_{y^{\vep}_{1}}$.
\eos
\bdim
For a given $v\in H^{1}$, a direct calculation yields
$$
|\nabla |v||^{2}=\frac12\left|\nabla v\right|^{2}+
\frac1{4|v|^{2}}\sum_{j=1}^N\left[(v_{j})^{2}(\bar v)^{2}+(v)^{2}(\bar v_{j})^{2}\right]=|\nabla v|^{2}-\frac{|\im(\bar v\nabla v)|^{2}}{|v|^{2}}
$$
where $v_{j}=v_{x_{j}}$ and, in the last term, it appears the square of the modulus of the vector whose 
components are $\im(\bar v v_{j})$. Then, we obtain
\begin{align*}
\E(\Upsilon^{\vep})=\E(|v^{\vep}_{1}|,|v^{\vep}_{2}|)+
\int\frac{|\im(\bar v^{\vep}_{1}\nabla v^{\vep}_{1})|^{2}}{|v^{\vep}_{1}|^{2}}dx
+\int\frac{|\im(\bar v^{\vep}_{2}\nabla v^{\vep}_{2})|^{2}}{|v^{\vep}_{2}|^{2}}dx.
\end{align*}
In turn, using Theorem~\ref{chiave}, it follows that, as $\eps$ vanishes,
\begin{align*}
0\leq \E(|v^{\vep}_{1}|,|v^{\vep}_{2}|) -&\E(R) 
+\int\frac{|\im(\bar v^{\vep}_{1}\nabla v^{\vep}_{1})|^{2}}{|v^{\vep}_{1}|^{2}}dx
+\int\frac{|\im(\bar v^{\vep}_{2}\nabla v^{\vep}_{2})|^{2}}{|v^{\vep}_{2}|^{2}}dx
\\
&\leq \alpha^{\vep}(t)+\eta^{\vep}(t)+\calo(\vep^{2}).
\end{align*}
Moreover, since $\|(|v^{\vep}_{1}|,|v^{\vep}_{2}|)\|_{2}=\|(v_{1}^{\vep},v_{2}^{\vep})\|_{2}=\|R\|_{2}$, we can conclude that
\be\label{dis1}
\int\frac{|\im(\bar v^{\vep}_{1}\nabla v^{\vep}_{1})|^{2}}
{|v^{\vep}_{1}|^{2}}dx+\int
\frac{|\im(\bar v^{\vep}_{2}\nabla v^{\vep}_{2})|^{2}}
{|v^{\vep}_{2}|^{2}}dx
\leq \alpha^{\vep}(t)+\eta^{\vep}(t)+\calo(\vep^{2}),
\ee
for every $\eps\in(0,\eps_{0})$ and for every $t\in [0,T^{\eps}_{*}]$.
Using  \eqref{veps} and \eqref{normel2}, for any $i=1,2$ we get
\begin{align*}
\frac{|\im(\bar v^{\vep}_{i}\nabla v^{\vep}_{i})|^{2}}
{|v^{\vep}_{i}|^{2}}
&=\frac{\left|\vep\im(\bar\phi^{\vep}_{i}(\vep x+x_{1},t)\nabla\phi^{\vep}_{i}(\vep x+x_{1},t))-\xi_{i}|\phi^{\vep}_{i}(\vep x+x_{1},t)|^{2}\right|^{2}}{|\phi^{\vep}_{i}(\vep x+x_{1},t)|^{2}}
\\
&=\vep^{2}\frac{\left|\im(\bar\phi^{\vep}_{i}(\vep x+x_{1},t)\nabla\phi^{\vep}_{i}(\vep x+x_{1}(t),t))\right|^{2}}{|\phi^{\vep}_{i}
(\vep x+x_{1},t)|^{2}}
+\xi_{i}^{2}|\phi^{\vep}_{i}(\vep x+x_{1},t)|^{2}
\\&\quad 
-2\vep\xi_{i} \im(\bar\phi^{\vep}_{i}(\vep x+x_{1},t)\nabla\phi^{\vep}_{i}(\vep x+x_{1}(t),t)).
\end{align*}
Whence, performing a change of variable and using definition \eqref{defpeps}, we derive
\be\label{equa1}
\int \frac{|\im(\bar v^{\vep}_{i}\nabla v^{\vep}_{i})|^{2}}
{|v^{\vep}_{i}|^{2}}dx
=\vep^{N}\int\frac{|p^{\vep}_{i}(x,t)|^{2}}{|\phi^{\vep}_{i}|^{2}}dx
+m_{i}\xi_i^2 -2 \xi_{i}\int p^{\vep}_{i}(x,t)dx.
\ee
Notice that
\begin{align*}
 \int\left|\vep^{N/2}\frac{p^{\vep}_{i}}{|\phi^{\vep}_{i}|}
-\frac{ \int p^{\vep}_{i}dx}{m_{i}}
\frac{|\phi^{\vep}_{i}|}{\vep^{N/2}}\right|^{2}dx
+m_{i}\left|\xi_{i}
-\frac{ \int p^{\vep}_{i}dx}{m_{i}}\right|^{2}
=\vep^{N}\int\frac{|p^{\vep}_{i}|^{2}}{|\phi^{\vep}_{i}|^{2}}dx
\\-\frac2{m_{i}}\left[\int p^{\vep}_{i}dx\right]^{2}+
\frac{\left[\int p^{\vep}_{i}\right]^{2}\int |\phi^{\vep}_{i}|^{2}}{\vep^{N}m_{i}^2}+m_{i}\xi_{i}^{2}
+\frac{\left[\int p^{\vep}_{i}\right]^{2}}{m_{i}}
-2\xi_{i}\int p^{\vep}_{i}dx
\end{align*}
which, by~\eqref{normel2} is equal to~\eqref{equa1}. In turn, \eqref{dis1} implies that
\begin{align}\label{dis2}
 \int\left|\vep^{N/2}\frac{p^{\vep}_{i}(x,t)}{|\phi^{\vep}_{i}|}
-\frac{ \int p^{\vep}_{i}dx}{m_{i}}
\frac{|\phi^{\vep}_{i}|}{\vep^{N/2}}\right|^{2}
&+m_{i}\left|\xi_{i}
-\frac{ \int p^{\vep}_{i}dx}{m_{i}}\right|^{2}
\\
\nonumber \leq \alpha^{\vep}(t)+\eta^{\vep}(t)&+\calo(\vep^{2}).
\end{align}
In order to prove the assertion, we need to estimate  $\rho^{\vep}_{i}(t)$ for $i=1,2$, where 
\beq\lbl{somma}
\rho^{\vep}_{i}(t)=\dys \left|\frac1{\vep^{N}}\int \psi(x)|\phi^{\vep}_{i}|^{2}dx-m_{i}\psi(y^{\vep}_{1})\right|+
\left|\int P^{\vep}(x,t)\psi(x)dx-M(t)\psi(y^{\vep}_{1})\right|
\eeq
for every function $\psi$ in $C^{2}$ such that 
$\|\psi\|_{C^{2}}\leq 1$. From the definition of \eqref{alfai} it holds
\begin{align*}
&\left|\int P^{\vep}(x,t)\psi(x)dx-M(t)\psi(y^{\vep}_{1})\right|\leq \\ 
&\leq \left|\int P^{\vep}(x,t)[\psi(x)-\psi(y^{\vep}_{1})]dx\right|+
\left| \psi(y^{\vep}_{1})\left(\int P^{\vep}(x,t)dx-M(t)\right)\right|
\\
&\leq \left|\int P^{\vep}(x,t)[\psi(x)-\psi(y^{\vep}_{1})]dx\right|
+|\alpha^{\vep}_{1}(t)+\alpha^{\vep}_{2}(t)|
\\
&\leq\sum_{i=1}^2\frac1{m_{i}}\left|\int p^{\vep}_{i}(x,t)dx\right|\left|\int\frac{\psi(x)|\phi^{\vep}_{i}(x,t)|^{2}}{\vep^{N}}dx-m_{i}\psi(y^{\vep}_{1})\right|
\\
&\;\,+\sum_{i=1}^2\left|\int \psi(x)\left[p^{\vep}_{i}(x,t)-\frac1{m_{i}}\left(\int p^{\vep}_{i}(x,t)dx\right) \frac{|\phi^{\vep}_{i}(x,t)|^{2}}{\vep^{N}}\right]dx
\right|
\\
&\;\,+|\alpha^{\vep}_{1}(t)+\alpha^{\vep}_{2}(t)|+\calo(\vep^{2}),
\end{align*}
for every $\eps\in(0,\eps_{0})$ and for every $t\in [0,T^{\eps}_{*}]$. Taking into account that $\int p^{\vep}_{i}dx$ is uni\-formly bounded and that, of course, 
$$
\int \Big[p^{\vep}_{i}(x,t)-\frac1{m_{i}}\Big(\int p^{\vep}_{i}(x,t)dx\Big)\frac{|\phi^{\vep}_{i}(x,t)|^2}{\vep^{N}}\Big]dx=0,
$$
there exists a positive constant $C_0$ such that, if we set 
$\tilde \psi(x)=\psi(x)-\psi(y^{\vep}_{1})$, it holds
\begin{align*}
\rho^{\vep}_{i}(t) &\leq 
\frac{1}{\vep^{N}}\int |\tilde\psi(x)||\phi^{\vep}_{i}(x,t)|^{2}dx+
\sum_{i=1}^2\frac{C_{0}}{\vep^{N}}\int |\tilde\psi(x)||\phi^{\vep}_{i}(x,t)|^{2}dx \\
&+\sum_{i=1}^2\int|\tilde\psi(x)|\left|p^{\vep}_{i}(x,t)
-\frac1{m_{i}}\Big(\int p^{\vep}_{i}(x,t)dx\Big)\frac{|\phi^{\vep}_{i}(x,t)|^{2}}{\vep^{N}}
\right|dx  \\
&+|\alpha^{\vep}_{1}(t)+\alpha^{\vep}_{2}(t)|+\calo(\vep^{2}).
\end{align*}
From Young inequality and \eqref{dis2} it follows (from now on $C_{0}$ 
will denote a constant that can vary from line to line)
\begin{align*}
 \rho^{\vep}_{i}(t) & \leq \frac{1}{\vep^{N}}\sum_{i=1}^2
 \int \left[ C_{0}|\tilde\psi(x)|  
+\frac12|\tilde \psi(x)|^{2}\right]
|\phi^{\vep}_{i}(x,t)|^{2}dx  \\
&+ \frac12\sum_{i=1}^2
\int \left|\frac{p^{\vep}_{i}(x,t)\vep^{N/2}}{|\phi^{\vep}_{i}(x,t)|}
-\frac1{m_{i}}\Big(\int p^{\vep}_{i}(x,t)dx\Big)\frac{|\phi^{\vep}_{i}(x,t)|}{\vep^{N/2}}
\right|^{2} \\
\noalign{\vskip4pt}
& +|\alpha^{\vep}_{1}(t)+\alpha^{\vep}_{2}(t)|+\calo(\vep^{2}) \\
&\leq  
\frac{1}{\vep^{N}}\sum_{i=1}^2\int
\left[ C_{0}|\tilde\psi(x)| +\frac12|\tilde \psi(x)|^{2}\right]
|\phi^{\vep}_{i}(x,t)|^{2} \\
&+|\alpha^{\vep}_{1}(t)+\alpha^{\vep}_{2}(t)|+\frac12\big[\alpha^{\vep}(t)+\eta^{\vep}(t)+\calo(\vep^{2})\big].
\end{align*}
Using the elementary inequality $a^{2}\leq 2b^{2}+2(a-b)^{2}$ 
with 
$$
a=\frac{\phi^{\vep}_{i}(x,t)}{\vep^{N/2}},\quad b=\frac1{\vep^{N/2}}
r_{i}\Big(\frac{x-y^{\vep}_{1}}{\vep}\Big),
$$
and recalling that $\tilde \psi$ is a uniformly bounded function
we derive
\begin{align*}
\rho^{\vep}_{i}(t) \leq & \frac1{\vep^{N}}\sum_{i=1}^2\int\left[ C_{0}|\tilde\psi(x)| +|\tilde \psi(x)|^{2}\right]
r_{i}^{2}\Big(\frac{x-y^{\vep}_{1}}{\vep}
\Big)dx \\
& +\frac{C_{0}}{\vep^{N}}\sum_{i=1}^2\int \Big|\phi^{\vep}_{i}(x,t)-
r_{i}\Big(\frac{x-y^{\vep}_{1}}{\vep}\Big)\Big|^{2}dx
\\
&+|\alpha^{\vep}_{1}(t)+\alpha^{\vep}_{2}(t)|+\frac12\big[\alpha^{\vep}(t)+\eta^{\vep}(t)+\calo(\vep^{2})\big].
\end{align*}
for every $\eps\in(0,\eps_{0})$ and for every $t\in [0,T^{\eps}_{*}]$.
Notice that $\tilde \psi$ satisfies the hypothesis of Lemma \ref{pote}
and $\tilde\psi (y^{\vep}_{1})=0$, then by virtue of inequality \eqref{quasifine}
we obtain the conclusion.
\edim
\vskip4pt

\subsection{Location estimates for $y^{\vep}_{1}$}
In the next results we start the study of the asymptotic behavior of 
$y^{\vep}_{1}$.

\begin{lemma}\lbl{propgamma}
Let us define the function
\beq\lbl{gamma}
\gamma^{\vep}(t)=\left|\gamma_{1}^{\vep}(t)\right|+\left|\gamma_{2}^{\vep}(t)\right|,
\;\; \text{with}\;\,
\gamma^{\vep}_{i}(t)=m_{i}x_{i}(t)-\frac1{\vep^{N}}
\int x\chi(x)|\phi^{\vep}_{i}(x,t)|^{2}dx,
\eeq
where $\chi(x)$ is the characteristic function defined in \eqref{defchi}.
Then $\gamma^{\vep}_{i}(t)$ is a continuous function with respect to
$t$ and  $\left|\gamma^{\vep}_{i}(0)\right|=\calo(\vep^{2})$ for $i=1,2$.
\end{lemma}
\bdim
The continuity of $\gamma^{\vep}$ immediately follows from the properties of the functions $\chi$ and $\phi^{\vep}_{i}$. In order to complete the proof,
note that Lemma \ref{pote} implies
$$
\quad\dys |\gamma^{\vep}_{i}(0)|=\dys \left|m_{i}\tilde x-\int (\tilde x+\vep y)\chi(\tilde x+\vep y)|r_{i}(y)|^{2}dy\right|
 \leq\dys C_{0}\vep^{2} +\left|m_{i}\tilde x-\tilde x\chi(\tilde x)\|r_{i}\|_{L^{2}}^{2}\right| 
$$
and as $\chi(\tilde x)=1$ we reach the conclusion.
\edim\vskip6pt

\begin{lemma}\lbl{y}
Let $T^{\vep}_{*}$ be the time introduced in \eqref{defteps}.
There exist positive constants $h_{0}$ and $\vep_{0}$ such that, if 
$|\eta_{i}^{\vep}|\leq h_{0}$ and $\vep\in(0,\vep_{0})$ 
there is a positive constant $L_{2}$ such that
$$
|x_{1}(t)-y^{\vep}_{1}(t)|\leq L_{2}
\left[\hat\alpha^{\vep}(t)+\eta^{\vep}(t)+\gamma^{\vep}(t)+\calo(\vep^{2})\right]
$$
for every $t\in [0,T^{\vep}_{*}]$.
\end{lemma}
\bdim
First we show that there exist $T_{0}>0$ and $B>0$ such that
\beq 
\lbl{stimay}
|y^{\vep}_{1}(t)|\leq B,
\eeq
for every $t\in [0,T^{\vep}_{*}]$ with 
$T^{\vep}_{*}\leq T_{0}$.  Let us first prove that
$$
\|\delta_{y^{\vep}_{1}(t_{2})}-\delta_{y^{\vep}_{1}(t_{1})}\|_{C^{2*}}<B,\quad\text{for all $t_{1},\,t_{2}\in [0,T^{\vep}_{*}]$}.
$$
Let $\psi\in C^{2}$ with $\|\psi\|_{C^2}\leq 1$ and pick $t_{1},\,t_{2}\in [0,T^{\vep}_{*}]$ 
with $t_{2}>t_{1}$.  From identity~\eqref{derfi} and integrating by parts, 
we obtain
\begin{align*}
\frac1{\vep^{N}}\int \psi(x)(|\phi^{\vep}_{i}(x,t_{2})|^{2}-|\phi^{\vep}_{i}(x,t_{1})|^{2})dx
&=\frac1{\vep^{N}}\int_{t_{1}}^{t_{2}}\!dt
\!\int \psi(x)\frac{\partial |\phi^{\vep}_{i}(x,t)|^2}{\partial t}dx
\\
&=-\int_{t_{1}}^{t_{2}}\!dt\!\int \psi(x){\rm div}\, p_{i}^{\vep}(x,t)dx
\\
&\leq \|\nabla \psi\|_{\infty}\int_{t_{1}}^{t_{2}}\!dt\!\int 
| p_{i}^{\vep}(x,t)|dx.
\end{align*}
It is readily seen from the $L^2$ estimate of $\nabla\phi^{\vep}_{i}$ 
that the last integral on the right hand side is uniformly bounded, so that
there exists a positive constant $C_{0}$ such that
$$
\||\phi^{\vep}_{i}(x,t_{2})|^2/\vep^{N}dx-|\phi^{\vep}_{i}(x,t_{1})|^2/\vep^{N}dx\|_{C^{2*}}\leq C_{0}|t_{2}-t_{1}|\leq C_{1}T_{0},
$$
with $C_{1}=2C_{0}$. Then Lemma~\ref{tecnico1},~\ref{propalfa},~\ref{propeta} 
and~\ref{propgamma} imply that the following inequality holds 
for sufficiently small $\vep$ and $h_0$ (the quantity $\alpha^\vep$ should be replaced
by $\hat\alpha^\vep$ in the definition of $T^{\vep}_*$)
\begin{align*}
m_{1}\|\delta_{y^{\vep}_{1}(t_{2})}-
\delta_{y^{\vep}_{1}(t_{1})}\|_{C^{2*}}
&\leq C_{1}T_{0}
+L[\hat\alpha^{\vep}(t_{2})+\eta^{\vep}(t_{2})+\hat\alpha^{\vep}(t_{1})
+\eta^{\vep}(t_{1})+\calo(\vep^{2})]
\\
&\leq C_{1}[T_{0}+ \calo(\vep^{2})+h_{0}].
\end{align*}
Here we fix $T_0$ and then $\vep_{0},\,h_{0}$ so small that $C_{1}\left[
T_{0}+ \calo(\vep^{2})+h_{0}\right]< M K_{0}$ where $K_{0}$ is 
the constant fixed in Lemma~\ref{propdelta} and from this lemma it follows
$$
|y^{\vep}_{1}(t_{2})-y^{\vep}_{1}(t_{1})|\leq  C_{2}K_{0},
$$
for every $t_{1},\,t_{2}\leq T_{0}$, and since $y^{\vep}_{1}(0)=\tilde x$
we obtain \eqref{stimay} for $B=C_{2}K_{0}+|\tilde x|$.
In view of property~\eqref{stimay} we can now prove the assertion.
Let us first observe that the properties of the function $\chi$ imply
\begin{align*}
|x_1(t)-y^{\vep}_{1}(t)|&=\frac1{m_1}\left|m_1x_1(t)-m_1
y^{\vep}_{1}(t)\right|
\\
&\leq 
\frac1{m_1}|\gamma^{\vep}_1(t)|
+\frac1{m_1}
\left|\int  \frac1{\vep^{N}} x\chi(x)|\phi^{\vep}_{1}(x,t)|^{2}-m_{1}y^{\vep}_{1}(t)\right|.
\end{align*}
Using \eqref{stimay} and \eqref{defchi} we obtain that
$\chi(y^{\vep}_{1})=1$, so that there exists a positive constant $C^{0}$ such that
\begin{align*}
|x_1(t)-y^{\vep}_{1}(t)| &\leq C_{0}\|x\chi\|_{C^{2}}
\||\phi^{\vep}_{1}|^{2}/{\vep^{N}}\,dx-m_{i}\delta_{y^{\vep}_{1}}\|_{C^{2*}}+C_{0}\gamma^{\vep}(t).
\end{align*}
This and Lemma \ref{tecnico1} give the conclusion.
\edim
\vskip6pt
In the previous Lemma we have fixed $T_{0}$ such that also Lemma \ref{tecnico1} and Theorem \ref{chiave} hold and 
now we are able to prove Theorem~\ref{mainth2}.
\vskip4pt
\noindent
{\bf Proof of Theorem \ref{mainth2}.}
We start the proof from the second conclusion of Theorem \ref{chiave}.
By Theorem \ref{chiave}, the family of continuous functions $\varrho^\vep:\R^+\to\R$, 
\begin{equation}
\label{Defrho}
\varrho^\vep(t)=\hat{L}\left[\hat\alpha^{\vep}(t)+\eta^{\vep}(t)+\gamma^{\vep}(t)\right]\quad \hat{L}=\max\{L,L_{1},L_{2}\}
\end{equation}
is such that $\varrho^\vep(0)=\calo(\vep^{2})$ and it satisfies 
\beq\lbl{fine}
{\Big \|}\Phi^{\vep}-\left(e^{\frac{{\rm i}}{\vep}(x\xi_{1}+\theta^{\vep}_{1})}
r_{1}\!\!\left(\frac{x-y^{\vep}_{1}}\vep\right),
e^{\frac{{\rm i}}{\vep}(x\xi_{2}+\theta^{\vep}_{2})}
r_{2}\!\!\left(\frac{x-y^{\vep}_{1}}\vep\right)\right)
{\Big \|}^{2}_{\H_{\vep}}\leq \varrho^\vep(t).
\eeq
Moreover, Lemma \ref{y} implies $|\vep \tilde y^{\vep}|=|x_{1}-y_{1}^{\vep}|\leq \varrho^\vep(t)$, so that $|\tilde y^{\vep}|=\frac{\varrho^\vep(t)}{\eps}$. Also,
$$
\|r_i(\cdot)-r_i(\cdot-\tilde y^{\vep})\|_{H^{1}}^{2}
\leq |\tilde y^{\vep}|^{2}\|\nabla r_i\|_{H^{1}}^{2}\leq C\left(\frac{\varrho^\vep(t)}{\eps}\right)^2,
\qquad \text{for all $i=1,2$}.
$$
Then
\begin{equation*}
{\Big \|}\Phi^{\vep}-\Big(e^{\frac{{\rm i}}{\vep}(x\cdot\xi_{1}+\theta^{\vep}_{1})}
r_1\Big(\frac{x-x_1(t)}\vep\Big),
e^{\frac{{\rm i}}{\vep}(x\cdot \xi_{2}+\theta^{\vep}_{2})}
r_2\Big(\frac{x-x_1(t)}\vep\Big)\Big)
{\Big \|}^{2}_{\H_{\vep}}\leq \tilde\varrho^\vep(t),
\end{equation*}
where we have set
$$
\tilde\varrho^\vep(t)=\varrho^\vep(t)+C\left(\varrho^\vep(t)/\eps\right)^2.
$$
Since $\tilde\varrho^\vep(0)=\calo(\vep^{2})$, the assertion follows.
\edim

\vskip2pt
\noindent
{\bf Proof of Theorem \ref{mainth2Bis}.}
In view of definition~\eqref{Defrho}, the assertion immediately 
follows by combining Lemmas \ref{tecnico1}, \ref{y} and \ref{propdelta}.\edim

\subsection{Smallness estimates for $\hat\alpha^\vep,\,\eta^\vep,\,\gamma^\vep$}
In the next lemma, under the assumptions of Theorem~\ref{mainth3}, we complete the study of
the asymptotic behaviour of system \eqref{Part-Sist} by obtaining the vanishing rate of
the functions $\hat\alpha^\vep$, $\eta^\vep$ and $\gamma^\vep$ as $\vep$ vanishes.
The time $T_0$ is the one chosen in the proof of Lemma~\ref{y}.

\begin{lemma}\lbl{limite}
Consider the framework of Theorem~\ref{mainth3}, that is $V=W$ and $\tilde\xi_1=\tilde\xi_2=\tilde\xi$.
Then there exists a positive constant $\bar L$ such that
$$
\hat\alpha^{\vep}(t)+\eta^{\vep}(t)+\gamma^{\vep}(t)\leq 
\bar L (T^{0})\vep^{2},
\quad \text{for every } t\in[0,T_{0}].
$$
\end{lemma}
\bdim
By the definition of $\alpha^\eps(t)$ (see formula~\eqref{defnalphaeta}) and taking into account that 
under the assumptions of Theorem~\ref{mainth3} it holds $\xi_1=\xi_2=\xi$ (with 
respect to the notations of Theorem~\ref{chiave}), there exists a positive constant
$C$ such that, for $t>0$,
$$
\hat\alpha^\eps(t)=\alpha^\vep(t)+|\alpha_1^{\vep}(t)+\alpha_2^{\vep}(t)|\leq (1+ |\xi(t)|)|\alpha_1^\eps(t)+\alpha_2^\eps(t)|\leq C |\alpha_1^\eps(t)+\alpha_2^\eps(t)|.
$$
Hence, without loss of generality, we can replace in the previous theorems (in 
particular Theorem~\ref{mainth3}) the quantities $\alpha^\vep(t)$ and $\hat\alpha^\eps(t)$ with
the absolute value $|\alpha_1^\eps(t)+\alpha_2^\eps(t)|$.
In a similar fashion, it is possible to replace the quantity $\gamma^\eps(t)$ defined in
formula \eqref{gamma} with the value $|\gamma_{1}^{\vep}(t)+\gamma_{2}^{\vep}(t)|$.
We will prove the desired assertion via Gronwall Lemma, so that we will first show that
there exists a positive constant $\bar L$ such that, for all $t\in[0,T^{\vep}_*]$,
\begin{align}
\lbl{alfalim}
\hat\alpha^{\vep}(t)\leq \calo(\vep^{2})+\bar L\int_{0}^{t}
[\hat\alpha^{\vep}(t)+\eta^{\vep}(t)+\gamma^{\vep}(t)]dt,
\\
\lbl{etalim}
\eta^{\vep}(t)\leq \calo(\vep^{2})+\bar L\int_{0}^{t}[\hat\alpha^{\vep}(t)+\eta^{\vep}(t)+\gamma^{\vep}(t)]dt,
\\
\lbl{gammalim}
\gamma^{\vep}(t)\leq \calo(\vep^{2})+\bar L\int_{0}^{t}
[\hat\alpha^{\vep}(t)+\eta^{\vep}(t)+\gamma^{\vep}(t)]dt.
\end{align}
Now, identity \eqref{derpot} of Lemma \ref{ideder} yield
\begin{align*}
\left|\frac{d}{dt}(\alpha^{\vep}_{1}+\alpha^{\vep}_{2})(t)\right|
& \leq \|\nabla V\|_{C^{2}}
\||\phi^{\vep}_{1}|^{2}/\vep^{N}-m_{1}\delta_{x(t)}\|_{C^{2*}} \\
&+\|\nabla V\|_{C^{2}}
\||\phi^{\vep}_{2}|^{2}/\vep^{N}-m_{2}\delta_{x(t)}\|_{C^{2*}}.
\end{align*}
Hence, using Lemmas \ref{tecnico1} and \ref{y} one obtains, for all $t\in[0,T^{\vep}_*]$,
\begin{align*}
\left|\frac{d}{dt}(\alpha^{\vep}_{1}+\alpha^{\vep}_{2})(t)\right|\leq 
L_{1}[\hat\alpha^{\vep}+\eta^{\vep}
+\gamma^{\vep}+\calo(\vep^{2})],
\end{align*}
for some positive constant $A_1$, yielding inequality~\eqref{alfalim}. As far as concern $\eta^{\vep}$, using \eqref{derfi} and Lemmas
\ref{tecnico1} and \ref{y} one has for  $t\in [0,T^{\vep}_*]$ that there exists
a positive constant $A_{2}$ such that, for all $t\in[0,T^{\vep}_*]$,
\begin{align*}
\left|\frac{d}{dt}(\eta^{\vep}_{1}+\eta^{\vep}_{2})(t)\right| &\leq
\Big|m_{1}\nabla V(x(t))\cdot \xi(t)+m_{2}\nabla V(x(t))\cdot\xi(t) \\
&+\int \chi(x)V(x){\rm div}_{x}p^{\vep}_{1}(x,t)+\int \chi(x)V(x){\rm div}_{x}p^{\vep}_{2}(x,t) \Big|
\\
&=\Big|\int \big[\nabla (\chi V)(p^{\vep}_{1}+p^\vep_{2})(x,t)-\nabla V(x(t))\cdot(m_1\xi(t)+m_2\xi(t))\big]dx\Big|
\\
\noalign{\vskip3pt}
&\leq \|\nabla (\chi V)\|_{C^{2}}\|P^{\vep}(x,t)dx-M(t)\delta_{x(t)}\|_{C^{2*}}
\\
\noalign{\vskip5pt}
&\leq A_{2}[\hat\alpha^{\vep}+\eta^{\vep}+\gamma^{\vep}+\calo(\vep^{2})]
\end{align*}
Let us now come to $\gamma^{\vep}$. By
the properties of the function $\chi$, identity \eqref{derfi}, Lemmas \ref{tecnico1} and \ref{y} it follows that 
there exists a positive constant $A_{3}$ such that, for all $t\in[0,T^{\vep}_*]$,
\begin{align*}
\Big|\frac{d}{dt}(\gamma^{\vep}_{1}+\gamma^{\vep}_{2})(t)\Big|&=\Big|\int \big[\nabla (x\chi)\cdot p^{\vep}_{1}(x,t)+\nabla (x\chi)\cdot p^{\vep}_{2}(x,t)\big]dx
-m_{1}\xi(t)-m_2\xi(t)\Big|\\
&=\Big|\int \big[\nabla (x\chi)\cdot P^{\vep}(x,t)-\nabla (x\chi)\cdot M(t)\delta_{x(t)}\big]dx\Big| \\
\noalign{\vskip4pt}
&\leq \|\nabla (x\chi)\|_{C^{2}}\|P^{\vep}(x,t)dx-M(t)\delta_{x(t)}\|_{C^{2*}} \\
\noalign{\vskip5pt}
&\leq A_{3}[\hat\alpha^{\vep}+\eta^{\vep}+\gamma^{\vep}
+\calo(\vep^{2})]
\end{align*}
Then inequalities \eqref{alfalim}, \eqref{etalim}, \eqref{gammalim} immediately
follow from Lemmas~\ref{propalfa}, \ref{propeta} and \ref{propgamma}.
The conclusion on $[0,T^\vep_*]$ is now a simple consequence of the Gronwall Lemma over $[0,T^\vep_*]$.
By the definition of $T^{\vep}_*$ and the continuity of $\alpha^\vep$, $\hat\alpha^\vep$ and $\eta^\vep$ we have that $T^{\vep}_*=T_0$ provided $\vep$ is chosen sufficiently small. 
To have this, one also has to take into account that, by construction 
(cf.\ formula~\eqref{Biggammaest}) and by the uniform smallness inequalities that we 
have just obtained over $[0,T^\vep_*]$, we reach
$$
\Gamma_{\Upsilon^\eps(t)}\leq 
L[\alpha^{\vep}(t)+\eta^{\vep}(t)+\calo(\vep^{2})]\leq \calo(\vep^{2}), \quad
\text{for all $t\in [0,T^\vep_*]$}.
$$
\edim

\section{Proofs of the main results}
\label{proofMR}

\noindent
{\bf Proof of Theorem \ref{mainth3}.}  In light of Lemma~\ref{limite} we have $\varrho^\eps(t),\tilde\varrho^\eps(t)\leq \bar L\eps^2$ for any $t\in [0,T_0]$.
Hence, the conclusions hold in $[0,T_{0}]$ as a direct consequence of 
Theorem \ref{mainth2}. Finally, taking as new initial data
$$
\phi^0_i(x):=r_i\big(\frac{x-x(T_0)}{\vep}\big) e^{\frac{{\rm i}}{\vep} x\cdot \xi(T_0)},
$$
and taking as a new a guiding Hamiltonian system
\begin{equation*}
\begin{cases}
\, \dot{\bar x}(t)=\bar\xi(t)
\\
\, \dot{\bar\xi}(t)=-\nabla V(\bar x(t))
\\
\bar x(0)=x(T_0),\,\,
\bar\xi(0)=\xi(T_0),
\end{cases}
\end{equation*}
the assertion is valid over $[T_0,2T_0]$. Reiterating ($T_0$ only depends on 
the problem) the argument yields the assertion
locally uniformly in time.
\edim

\noindent
{\bf Proof of Theorem \ref{mainth-Bis}.} 
Combining definition~\eqref{Defrho} with the assertions 
of Lemmas~\ref{tecnico1} and~\ref{limite}, we obtain the property 
over the interval $[0,T_{0}]$. Then we can argue as in the proof of 
Theorem~\ref{mainth3} to achieve the conclusion locally uniformly in time. \edim

\bigskip

\bigskip

\end{document}